\documentclass{article}

\usepackage{microtype}
\usepackage{graphicx}
\usepackage{subfigure}
\usepackage{booktabs}
\usepackage{amsmath}
\usepackage{graphics}
\usepackage{color}
\usepackage{mathrsfs}
\usepackage{amssymb}
\usepackage{amsmath}
\usepackage{boxedminipage}
\usepackage{stmaryrd}
\usepackage{stfloats}
\usepackage{multirow}
\usepackage{accents} 
\newtheorem{definition}{Definition}
\newtheorem{proposition}{Proposition}

\newtheorem{lemma}{Lemma}
\newtheorem{theorem}{Theorem}
\newtheorem{corollary}{Corollary}
\newtheorem{assumption}{Assumption}
\newtheorem{algorithm}{Algorithm}

\usepackage{hyperref}
\def\proof{\par\noindent{\em Proof.}}
\def\endproof{\hfill $\Box$ \vskip 0.4cm}
\newcommand{\RR}{\mathbf R}

\usepackage[accepted]{icml2020}

\icmltitlerunning{Linear Convergence of RPDC method for Large-scale LCCP}

\begin{document}

\twocolumn[
\icmltitle{Linear Convergence of Randomized Primal-Dual Coordinate Method for Large-scale Linear Constrained Convex Programming}



\icmlsetsymbol{equal}{*}

\begin{icmlauthorlist}
\icmlauthor{Daoli Zhu}{equal,at}
\icmlauthor{Lei Zhao}{equal,no}
\end{icmlauthorlist}

\icmlaffiliation{at}{Antai College of Economics and Management and Sino-US Global Logistics
Institute, Shanghai Jiao Tong University, Shanghai, China}
\icmlaffiliation{no}{School of Naval Architecture, Ocean and Civil Engineering, Shanghai
Jiao Tong University, Shanghai, China}

\icmlcorrespondingauthor{Lei Zhao}{l.zhao@sjtu.edu.cn}

\icmlkeywords{Linear constrained convex programming, Support vector machine, Machine learning portfolio, Randomized coordinate primal-dual method, Linear convergence}

\vskip 0.3in
]



\printAffiliationsAndNotice{\icmlEqualContribution} 

\begin{abstract}
Linear constrained convex programming has many practical applications, including support vector machine and machine learning portfolio problems. We propose the randomized primal-dual coordinate (RPDC) method, a randomized coordinate extension of the first-order primal-dual method by~\cite{CohenZ} and~\cite{ZhaoZhu2017}, to solve linear constrained convex programming. We randomly choose a block of variables based on a uniform distribution, linearize, and apply a Bregman-like function (core function) to the selected block to obtain simple parallel primal-dual decomposition. We then establish almost surely convergence and expected $O(1/t)$ convergence rate, and expected linear convergence under global strong metric subregularity. Finally, we discuss implementation details for the randomized primal-dual coordinate approach and present numerical experiments on support vector machine and machine learning portfolio problems to verify the linear convergence.
\end{abstract}
\section{Introduction}
\label{intro}
This paper considers linear constrained convex programming (LCCP),
\begin{equation}\label{Prob:general-function}
\begin{array}{lll}
\mbox{(P):}  &\min       & F(u)=G(u)+J(u)      \\
             &\rm {s.t}  & Au-b=0 \\
             &           & u\in \mathbf{U}
\end{array},
\end{equation}
where $G$ is a convex smooth function on the closed convex set $\mathbf{U}\subset \RR^{n}$; and $J$ is a convex, possibly non-smooth function on $\mathbf{U}\subset \RR^{n}$. We assume that $J(u)=\sum\limits_{i=1}\limits^{N}J_i(u_i)$ is additive with respect to the space decomposition
\begin{equation}\label{spacedecomposition}
\mathbf{U}=\mathbf{U}_1\times\mathbf{U}_2\cdots\times\mathbf{U}_N, u_i\in\mathbf{U}_i\subset \RR^{n_i}~\mbox{and}~ \sum\limits_{i=1}\limits^{N}n_i=n.
\end{equation}
Each $J_{i}$ is a convex but possibly non-smooth function on $U_{i}\subset \RR^{n_{i}}$; $A=(A_{1},A_{2},\cdots,A_{N})\in \RR^{m\times n}$ is an appropriate partition of $A$, where $A_{i}$ is an $m\times n_{i}$ matrix and $b\in \RR^{m}$ is a vector.
\subsection{Motivation}
Linear constrained convex programming is an important and challenging application problem class. We present several example applications to demonstrate the reasons for interest in type (P) problems.
\subsubsection{Support vector machine}
Support vector machine (SVM) is a popular supervised learning method~\cite{boser1992training,cortes1995support}, widely used for pattern recognition~\cite{burges1998tutorial,scholkopf2000new} and classification~\cite{chang2011libsvm}. The SVM problem can be expressed as
\begin{eqnarray*}
\begin{array}{lcl}
\mbox{\rm(SVM)}&\min\limits_{u\in[0,c]^n}   & \frac{1}{2} u^{\top}Q u-\mathbf{1}_n^{\top}u  \\
               &\mbox{s.t.}                 &  y^{\top}u=0
\end{array},
\end{eqnarray*}
where $u\in\RR^{n}$ are the decision variables, $Q\in\RR^{n\times n}$ is a symmetric and positive-definite matrix, $c\in\RR$ is the upperbound of all variables, $y\in\{-1,1\}^n$ is the vector of labels, and $\mathbf{1}_n$ is an $n$-dimensional vector of $1$s.
\subsubsection{Machine learning portfolio problem}
Portfolio optimization (PO) via machine learning has received increased attention recently. PO aims to invest in a group of financial assets with instructions by machine, based on financial principles and optimization strategies. Since this requires considerable quantitative calculation, machine learning methods are essential to reduce human mistakes and biases for real-world investment.~\cite{portfolio1,portfolio2,portfolio3,portfolio5,portfolio4} The machine learning portfolio (MLP) problem can be expressed as
\begin{eqnarray*}
\begin{array}{lcl}
\mbox{\rm(MLP)}&\min\limits_{u\in\RR^n}   & \frac{1}{2} u^{\top}\Sigma u + \lambda\|u\|_1 \\
               &\mbox{s.t.}                 &  \mu^{\top}u=\rho \\
               &                            &  \mathbf{1}_n^{\top}u=1 \\
\end{array},
\end{eqnarray*}
where $u\in \RR^{n}$ is the decision portfolio vector; $\Sigma\in \RR^{n\times n}$, which is symmetric and positive-definite, is the estimated covariance matrix of asset returns; $\mu\in\RR^{n}$ is the expectation of asset returns; $\rho$ is a predefined prospective growth rate; and $\mathbf{1}_n$ is an $n$-dimensional vector of $1$s.
\subsection{Related works}
The augmented Lagrangian method (ALM) is an approach for general LCCP that can overcome the dual Lagrangian instability and non-differentiability. The augmented Lagrangian for a constrained convex program has the same solution set as the original constrained convex program. Consider the following augmented Lagrangian function of (P),
\begin{equation}\label{func:AL0}
L_{\gamma}(u,p)=F(u)+\langle p,Au-b\rangle+\frac{\gamma}{2}\|Au-b\|^{2}.
\end{equation}
The ALM for equality-constrained optimization problems can be expressed as~\cite{Hestenes1969,Powell1969}
\begin{eqnarray*}
\begin{array}{l}
\mbox{{\bf ALM}}\\
\left\{\begin{array}{l}
u^{k+1}=\arg\min\limits_{u\in \mathbf{U}}L_{\gamma}(u,p^{k});\label{primal_ALM}\\
p^{k+1}= p^k+\gamma(Au^{k+1}-b).\label{dual_ALM}
\end{array}
\right.
\end{array}
\end{eqnarray*}
Although ALM has several advantages, it does not preserve separability, even when the initial problem is separable. One way to decompose the augmented Lagrangian is to use alternating direction method of multipliers (ADMM)~\cite{ADMM1983}, which applies a Gauss-Seidel-like minimization strategy. Another way to overcome this difficulty is the auxiliary problem principle of augmented Lagrangian (APP-AL)~\cite{CohenZ} method, a fairly general first-order primal-dual decomposition method based on linearizing the augmented Lagrangian of LCCP. The APP-AL scheme for LCCP can be expressed as
\begin{eqnarray*}
\begin{array}{l}
\mbox{{\bf APP-AL}}\\
\left\{\begin{array}{l}
u^{k+1}=\arg\min\limits_{u\in \mathbf{U}}\langle\nabla G(u^{k}), u\rangle + J(u)+ \langle q^k,Au\rangle\nonumber\\
\qquad\qquad\qquad\qquad\qquad\qquad\qquad\quad+\frac{1}{\epsilon}D(u,u^k);\label{primal_APP}\\
p^{k+1}= p^k+\gamma(Au^{k+1}-b),\label{dual_APP}
\end{array}
\right.
\end{array}
\end{eqnarray*}
\indent where $q^k=p^k+\gamma(Au^k-b)$, $D(u,v)=K(u)-K(v)-\langle\nabla K(v),u-v\rangle$ is a Bregman like function with $K$ is strongly convex and gradient Lipschitz.~\cite{ZhaoZhu2017} extended~\cite{CohenZ} to propose first-order primal-dual augmented Lagrangian methods for nonlinear convex cone programming with separable and non-separable objective and constraints as an algorithm (variant auxiliary problem principle, VAPP). They showed that APP-AL can be viewed as a forward-backward splitting method to find a solution for (P): $w^{k+1}=(\Gamma^k+\mathcal{A})^{-1}(\Gamma^k-\mathcal{B})w^k$, where
\begin{equation*}
\mathcal{A}(w)=\left(\begin{array}{c}\partial J(u)+\mathcal{N}_{\mathbf{U}}(u)\\ \mathbf{0}_m\end{array}\right)
\end{equation*}
and
\begin{equation*}
\mathcal{B}(w)=\left(\begin{array}{c}\nabla G(u)+A^{\top}[p+\gamma(Au-b)]\\ b-Au\end{array}\right)
\end{equation*}
are two maximal monotone mappings; $w=(u,p)$; $\mathcal{N}_{\mathbf{U}}(u)=\{\xi:\langle\xi,\zeta-u\rangle\leq0,\forall\zeta\in\mathbf{U}\}$ is the normal cone at $u$ to $\mathbf{U}$; and $\Gamma^k(w)=\left(\begin{array}{c}\frac{1}{\epsilon}\nabla K(u)\\ \frac{1}{\gamma}[p-p^k]-[Au-b]\end{array}\right)$ is a strongly monotone nonlinear Bregman operator. Problem (P) can be reformulated as an augmented Lagrangian based inclusion problem,
\begin{equation}\label{FBS_inclusion}
\begin{array}{l}
0\in H_{\gamma}(w)=\left(\begin{array}{c}
\partial_uL_{\gamma}(u,p)+\mathcal{N}_{\mathbf{U}}(u)
\\
\nabla_pL_{\gamma}(u,p)
\end{array}\right)\\
=\left(\begin{array}{c}
\nabla G(u)+\partial J(u)+A^{\top}[p+\gamma(Au-b)]\\
\qquad\qquad\qquad\qquad\qquad\qquad+\mathcal{N}_{\mathbf{U}}(u)
\\
b-Au
\end{array}\right)\\
=\mathcal{A}(w)+\mathcal{B}(w).
\end{array}
\end{equation}
We propose an iteration based forward-backward splitting algorithm to solve~\eqref{FBS_inclusion}, $w^{k+1}=(\Gamma^k+\mathcal{A})^{-1}(\Gamma^k-\mathcal{B})w^k$, hence $(\Gamma^k-\mathcal{B})w^k\in(\Gamma^k+\mathcal{A})w^{k+1}$. Thus,
\begin{equation*}
0\in\left(\begin{array}{l}\frac{1}{\epsilon}[\nabla K(u^{k+1})-\nabla K(u^k)]+\nabla G(u^k)\\
\qquad+A^{\top}[p^k+\gamma(Au^k-b)]+\partial J(u^{k+1})\\
\qquad\qquad\qquad\qquad\qquad\qquad+\mathcal{N}_{\mathbf{U}}(u^{k+1})\\
p^{k+1}-p^k-\gamma(Au^{k+1}-b)\end{array}\right),
\end{equation*}
which is the exact inclusion problem formulation for APP-AL. Zhao and Zhu also proposed linear convergence for VAPP~\cite{ZhaoZhu2017}.

Recent big data applications employ very large datasets that are commonly distributed over different locations. Hence it is often impractical to assume that optimization algorithms can traverse an entire dataset once in each iteration, because this would be very time consuming and/or unreliable, often resulting in low resource utilization due to synchronization among different computing units, e.g. CPUs, GPUs, and cores, in a distributed computing environment. On the other hand, block coordinate descent (BCD) algorithms can make progress by using information obtained from a randomly selected data subsets and, hence provide more flexible implementation in distributed environments. The main advantage of BCD is to reduce complexity and memory requirements per iteration, which becomes increasingly important for very-large scale problems.

We briefly review some related works on coordinate type methods.

Two BCD variations are widely employed for problems without constraints. The first BCD variation relates to block-choosing strategy. One common approach is a cyclic strategy.~\cite{Tseng} proved BCD cyclic strategy convergence and~\cite{LuoZeng} and~\cite{WangLin} proved local and global linear convergence, respectively, under specific assumptions. Another common approach is a randomized strategy.~\cite{Nesterov} studied the convergence rate of randomized BCD for convex smooth optimization and~\cite{RT} and~\cite{LuXiao} subsequently extended Nesterov's technique to composite optimization. The second BCD variation relates to the point read to evaluate the gradient in each iteration. The approaches are called asynchronous BCD if the read points have different "ages", and synchronous BCD otherwise. All BCD variants reviewed above are synchronous BCD.~\cite{LiuWright1} and~\cite{LiuWright2} established the convergence rate of asynchronous BCD for composite optimization and convex smooth optimization, respectively, without constraints.

Few previous studies considered BCD methods for problems with constraints.~\cite{Necoara} proposed a random coordinate descent algorithm for an optimization problem with one linear constraint. ~\cite{GaoXuZhang} and~\cite{Xu18} considered a similar scheme to RPDC, obtaining expected $O(1/t)$ and $O(1/t^2)$ rates, respectively.~\cite{Xu19} recently proposed an asynchronous RPDC algorithm with expected $O(1/t)$ rate. However, to the best of our knowledge, no previous study considered convergence and linear convergence results for RPDC.

This paper focused on RPDC, the randomized coordinate extension of APP-AL, as shown in Algorithm~\ref{alg:RPDC}.
\subsection{Main contributions and outline for this paper}
\indent We propose the randomized primal-dual coordinate (RPDC) method based on the first-order primal-dual method~\cite{CohenZ,ZhaoZhu2017}. The RPDC method randomly updates one block of variables based on a uniform distribution. The main contributions from this paper are as follows:
\begin{itemize}
\item[{\rm(i)}] We show that the sequence generated by RPDC converges to an optimal solution with probability $1$.
\item[{\rm(ii)}] We show RPDC has expected $O(1/t)$ rate for general LCCP.
\item[{\rm(iii)}] We establish the expected linear convergence of RPDC under global strong metric subregularity.
\item[{\rm(iv)}] We show that SVM and MLP problems satisfy global strong metric subregularity under some reasonable conditions.
\item[{\rm(v)}] Finally, we discuss the implementation details of RPDC and present numerical experiments on SVM and MLP problems to verify linear convergence.
\end{itemize}
The remainder of this paper is organized as follows. Section~\ref{Pre} discusses technical preliminaries. Section~\ref{convergence_RPDC} shows almost surely convergence and expected $O(1/t)$ convergence rate for RPDC. Section~\ref{rate} establishes the expected linear convergence of RPDC under global strong metric subregularity. Section~\ref{num} discusses implementation details for RPDC and presents numerical experiments on SVM and MLP problems. Finally, Section~\ref{ccl} summarizes and concludes the paper.
\section{Preliminaries}\label{Pre}
This section provides some useful preliminaries for subsequent discussions and summarizes notations and assumptions. We denote vector inner product and Euclidean norm as $\langle\cdot\rangle$ and $\|\cdot\|$, respectively.
\subsection{Notations and assumptions}
Throughout this paper, we make the following standard assumptions for Problem (P).
\begin{assumption}\label{assump1}
{\rm(H$_1$)} $J$ is a convex, lower semi-continuous function (not necessarily differentiable) such that $\mathbf{dom}J\cap \mathbf{U}\neq \emptyset$.\\
{\rm(H$_2$)} $G$ is convex and differentiable, and its derivative is Lipschitz with constant $B_{G}$.\\
{\rm(H$_3$)} There exists at least one saddle point for the Lagrangian of {\rm(P)}.
\end{assumption}
From Assumption~\ref{assump1} and Theorem 3.2.12~\cite{Ortega70}, the following descent property for $G$ holds
\begin{eqnarray}
G(v)-G(u)-\langle\nabla G(u),v-u\rangle\leq\frac{B_G}{2}\|u-v\|^2.\label{eq:Lip_G}
\end{eqnarray}
\subsection{Lagrangian and Karush-Kuhn-Tucker mapping}
The Lagrangian of (P) is defined as
\begin{equation}\label{func:L}
L(u,p)=F(u)+\langle p,Au-b\rangle,
\end{equation}
and a saddle point $(u^*,p^*)\in \mathbf{U}\times\mathbf{R}^m$ is such that
\begin{equation}
\forall u\in\mathbf{U}, p\in \mathbf{R}^{m}: L(u^{*},p)\leq L(u^{*},p^{*})\leq L(u,p^{*}). \label{saddle point:L}
\end{equation}
From Assumption~\ref{assump1}, there exist saddle points of $L$ on $\mathbf{U}\times\mathbf{R}^m$, and we denote the set of saddle points as $\mathbf{U}^*\times\mathbf{P}^*$. By definition, saddle point $(u,p)\in\mathbf{U}^*\times\mathbf{P}^*$ of $L$ satisfies
\begin{equation}\label{eq:L_inclusion}
\left\{
\begin{array}{l}
0\in\partial_uL(u,p)+\mathcal{N}_{\mathbf{U}}(u);\\
0=-\nabla_pL(u,p).
\end{array}
\right.
\end{equation}
System~\eqref{eq:L_inclusion} can also be considered the Karush-Kuhn-Tucker (KKT) system of (P). Thus, the saddle point problem of (P) can be represented as the inclusion problem
\begin{eqnarray*}
0\in H(w)&=&\left(
\begin{array}{c}
\partial_uL(u,p)+\mathcal{N}_{\mathbf{U}}(u) \\
-\nabla_pL(u,p)
\end{array}
\right)\\
&=&\left(
\begin{array}{c}
\nabla G(u)+\partial J(u)+A^{\top}p+\mathcal{N}_{\mathbf{U}}(u) \\
b-Au
\end{array}
\right),
\end{eqnarray*}
where, we call $H$ the KKT mapping for obvious reasons.
\section{Convergence and convergence rate analysis of RPDC}\label{convergence_RPDC}
This section establishes almost surely convergence and expected $O(1/t)$ convergence rate for RPDC. First, we introduce the following assumption on core function $K$ and parameters $\epsilon$ and $\rho$:
\begin{assumption}\label{assump2}
\begin{itemize}
\item[{\rm(i)}] $K$ is strongly convex with parameter $\beta$ and differentiable with its gradient Lipschitz continuous with parameter $B$ on $\mathbf{U}$.
\item[{\rm(ii)}] Parameters $\epsilon$ and $\rho$ satisfy:
\begin{equation}\label{para-choice-convergence}
0<\epsilon<\beta/[B_G+\gamma\lambda_{\max}(A^{\top}A)]\mbox{and}\;0<\rho<\frac{2\gamma}{2N-1},
\end{equation}
where $\lambda_{\max}(A^{\top}A)$ is the largest eigenvalue of $A^{\top}A$.
\end{itemize}
\end{assumption}
\indent Let $D(u,v)=K(u)-K(v)-\langle\nabla K(v),u-v\rangle$ be a Bregman like function (core function)~\cite{Mirror,CohenZ}. Two popular core functions $K$ satisfy Assumption~\ref{assump2}:
\begin{enumerate}
\item $K(u)=\frac{1}{2}\|u\|^2$, where $\beta=B=1$, and
\item $K(u)=\frac{1}{2}\|u\|_{Q}^2$, where $Q$ is the $Q$-quadratic norm associated with positive definite matrix $Q$.
\end{enumerate}
From Assumption~\ref{assump2}, $\frac{\beta}{2}\|u-v\|^2\leq D(u,v)\leq\frac{B}{2}\|u-v\|^2$ and Algorithm~\ref{alg:RPDC} shows the proposed RPDC method to solve (P). For the sake of brevity, let us set $q^k=p^k+\gamma(Au^k-b)$. Then the primal problem can be expressed as
\begin{equation*}
\begin{array}{l}
\mbox{\rm (AP$^k$)}\;\min\limits_{u\in\mathbf{U}}\langle\nabla_{i(k)}G(u^k),u_{i(k)}\rangle+J_{i(k)}(u_{i(k)})\\
\qquad\qquad\quad+\langle q^k,A_{i(k)}u_{i(k)}\rangle+\frac{1}{\epsilon}[K(u)-\langle\nabla K(u^k),u\rangle].
\end{array}
\end{equation*}
If we choose an additive Bregman function (or core function) with respect to the space decomposition~\eqref{spacedecomposition}, i.e., $K(u)=\sum\limits_{i=1}^{N}K_i(u_i)$, then problem (AP$^k$) is just a small optimization problem for selected block $i(k)$. Thus, taking $K(u)=\sum\limits_{i=1}^N\frac{\|u\|^2}{2}$ for (AP$^k$), we perform only a block proximal gradient update for block $i(k)$, where we linearize the coupled function $G(u)$ and augmented Lagrangian term $\langle p,Au-b\rangle+\frac{\gamma}{2}\|Au-b\|^2$, and add the proximal term to it.\\
\begin{algorithm}[tb]
   \caption{Proposed randomized primal-dual coordinate method}
   \label{alg:RPDC}
\begin{algorithmic}
   \FOR{$k=1$ {\bfseries to} $t$}
   \STATE Choose $i(k)$ from $\{1,\ldots,N\}$ with equal probability;
   \STATE $u^{k+1}=\arg\min\limits_{u\in \mathbf{U}}\langle\nabla_{i(k)} G(u^{k}), u_{i(k)} \rangle + J_{i(k)}(u_{i(k)})$
   \STATE \qquad\qquad\qquad\qquad$+\langle q^k,A_{i(k)}u_{i(k)}\rangle+\frac{1}{\epsilon}D(u,u^k)$;
   \STATE $p^{k+1}= p^{k}+\rho(Au^{k+1}-b)$.
   \ENDFOR
\end{algorithmic}
\end{algorithm}
Indices $i(k)$, $k=0,1,2,\ldots$ in Algorithm~\ref{alg:RPDC} are random variables. After $k$ iterations, RPDC generates random output $(u^{k+1}, p^{k+1})$. We denote $\mathcal{F}_k$ as a filtration generated by the random variable $i(0),i(1),\ldots,i(k)$, i.e.,
$$\mathcal{F}_{k}\overset{def}{=}\{i(0),i(1),\ldots,i(k)\}, \mathcal{F}_{k}\subset\mathcal{F}_{k+1},$$
and define $\mathcal{F}=(\mathcal{F}_{k})_{k\in\mathbb{N}}$,  $\mathbb{E}_{\mathcal{F}_{k+1}}=\mathbb{E}(\cdot|\mathcal{F}_{k})$ as the conditional expectation with respect to $\mathcal{F}_{k}$. The conditional expectation in the $i(k)$ term for given $i(0),i(1),\ldots,i(k-1)$ is $\mathbb{E}_{i(k)}$.

Let $w=(u,p)$, given $w^*=(u^*,p^*)\in\mathbf{U}^*\times\mathbf{P}^*$. Then for any $w,w'\in\mathbf{U}\times\RR^m$, we construct the function
\begin{eqnarray}\label{Lambda}
\Lambda(w,w')&=&D(u',u)+\frac{\epsilon}{2N\rho}\|p-p'\|^2\nonumber\\
&&+\frac{\epsilon(N-1)}{N}[L(u,p)-L(u^*,p^*)]\nonumber\\
&&+\frac{\epsilon(N-2)\gamma}{2N}\|Au-b\|^2
\end{eqnarray}
and we have the following lemma regarding the boundness of $\Lambda(w,w^*)$ and $\Lambda(w,w')$.
\begin{lemma}[Boundness of $\Lambda(w,w^*)$ and $\Lambda(w,w')$]\label{boundness_Lambda}
Suppose Assumption \ref{assump1} and \ref{assump2} hold. $w^*=(u^*,p^*)\in\mathbf{U}^*\times\mathbf{P}^*$. Then there exist positive numbers $d_1$, $d_2$ and $d_3$, such that
\begin{itemize}
\item[{\rm(i)}] $\Lambda(w,w^*)\geq d_1\|w-w^*\|^2$,
\item[{\rm(ii)}] $\Lambda(w,w^*)\leq d_2\|w-w^*\|^2
+\frac{\epsilon(N-1)}{N}[L(u,p^*)-L(u^*,p^*)]$, and
\item[{\rm(iii)}] $\Lambda(w,w')\geq-d_3\|p-p^*\|^2$;
\end{itemize}
with $d_1=\min\bigg{\{}\frac{1}{2N}[N\beta-\epsilon\gamma\lambda_{\max}(A^{\top}A)],\frac{\epsilon}{4N\gamma}\bigg{\}}$,\\ $d_2=\max\bigg{\{}\frac{(4N-3)\epsilon}{(4N-2)N\rho},\frac{NB+\epsilon(2N-3)\gamma\lambda_{\max}(A^{\top}A)}{2N}\bigg{\}}$,\\
 and $d_3=\frac{\epsilon(N-1)^2}{2\gamma N(N-2)}$.
\end{lemma}
To analyze the convergence of RPDC, we need the point $T(w^k)=\big{(}T_u(w^k),T_p(w^k)\big{)}$ generated by one deterministic iteration of APP-AL for given $w^k$,
\begin{equation*}
\begin{array}{l}
\mbox{{\bf APP-AL}}\\
\left\{
\begin{array}{l}
T_u(w^k)=\arg\min\limits_{u\in \mathbf{U}}\langle\nabla G(u^{k}), u\rangle+ J(u)+\langle q^k,Au\rangle\\
\qquad\qquad\qquad\qquad\qquad\qquad\qquad\qquad+\frac{1}{\epsilon}D(u,u^k);\\
T_p(w^k)= p^k+\gamma\left[AT_u(w^k)-b\right],
\end{array}
\right.
\end{array}
\end{equation*}
and the following lemma.
\begin{lemma}\label{lemma:bound1} {\bf (Estimation on the variance of $\Lambda(w^k,w)$)}
Let Assumption~\ref{assump1} and~\ref{assump2} hold and $\{(u^k,p^k)\}$ be generated by RPDC. Then
\begin{eqnarray*}
&&\Lambda(w^k,w)-\mathbb{E}_{i(k)}\Lambda(w^{k+1},w)\nonumber\\
&\geq&\frac{\epsilon}{N}\mathbb{E}_{i(k)}\big{[}L(u^{k+1},p)-L(u,q^{k})\big{]}+d_4\|w^k-T(w^k)\|^2,
\end{eqnarray*}
where $d_4=\frac{\min\bigg{\{}\frac{\beta-\epsilon[B_G+\gamma\lambda_{\max}(A^{\top}A)]}{2},\frac{\epsilon[2\gamma-(2N-1)\rho]}{2N}\bigg{\}}}{\max\{N^2+2\gamma^2(N^2+2)\lambda_{\max}(A^{\top}A),4\gamma^2\}}$.
\end{lemma}
\subsection{Almost surely convergence of RPDC}
Based Lemma~\ref{lemma:bound1}, we establish the convergence of RPDC.
\begin{theorem}[Almost surely convergence]\label{theo:convergence}
Let the assumptions of Lemma~\ref{lemma:bound1} hold, then
\begin{itemize}
\item[{\rm(i)}] $\sum\limits_{k=0}\limits^{+\infty}\|w^k-T(w^k)\|^2<+\infty$ $\mbox{a.s.}$.
\item[{\rm(ii)}] The sequence $\{w^{k}\}$ generated by RPDC is almost surely bounded.
\item[{\rm(iii)}] Every cluster point of $\{w^{k}\}$ almost surely is a saddle point of the Lagrangian for (P).
\end{itemize}
\end{theorem}
\subsection{Convergence rate analysis for RPDC}
This subsection provides the convergence rate of RPDC. We define the average sequence for the sequence $\{(u^k,p^k)\}$ generated from Algorithm RPDC and any $t>0$ as
$$\bar{u}_{t}=\frac{\sum_{k=0}^{t}u^{k+1}}{t+1}\;\mbox{and}\;\bar{p}_{t}=\frac{\sum_{k=0}^{t}q^k}{t+1}.$$
\begin{theorem}\label{thm:primal_rate} {\bf(Expected primal suboptimality and expected feasibility)}\\
\indent Let Assumption~\ref{assump1} and~\ref{assump2} hold; $(u^*,p^*)$ be a saddle point, and $M$ be a bound of dual optimal solution of (P), i.e., $\|p^*\|<M$; and $\{(u^k,p^k)\}$ be generated by RPDC. Then we have the following identities.
\begin{itemize}
\item[{\rm(i)}]Global estimate of expected bifunction values.
$$
\mathbb{E}_{\mathcal{F}_t}\big{[}L(\bar{u}_t,p)-L(u,\bar{p}_t)\big{]}\leq\frac{Nh(w^0,w)}{\epsilon(t+1)},$$
$\forall u\in\mathbf{U}, p\in\RR^m$, $(u,p)$ could possibly be random and $h(w,w')=\Lambda(w,w')+\frac{d_3}{d_1}\Lambda(w,w^*)\geq0$.
\item[{\rm(ii)}] Expected feasibility.
$$\mathbb{E}_{\mathcal{F}_{t}}\|A\bar{u}_t-b\|\leq\frac{Nd_5}{(M-\|p^*\|)\epsilon(t+1)},$$
where $d_5=\sup\limits_{\|p\|<M}h(w^0,(u^*,p))$.
\item[{\rm(iii)}] Expected suboptimality.
\begin{eqnarray*}
\mathbb{E}_{\mathcal{F}_{t}}\left[F(\bar{u}_t)-F(u^*)\right]&\geq&-\frac{\|p^*\|Nd_5}{(M-\|p^*\|)\epsilon(t+1)},\\
\mathbb{E}_{\mathcal{F}_{t}}\left[F(\bar{u}_t)-F(u^*)\right]&\leq&\frac{Nd_5}{\epsilon(t+1)}.
\end{eqnarray*}
\end{itemize}
\end{theorem}
\section{Linear convergence of RPDC under global strong metric subregularity}\label{rate}
This section establishes the expected linear convergence of RPDC. Let $\phi(w,w^*)=\Lambda(w,w^*)+\frac{\epsilon}{N}[L(u,p^*)-L(u^*,p^*)]$. Then $\phi(w,w^*)\geq 0$ and $\phi(w,w^*)=0$ if and only if $w=w^*$, and following lemma holds for the upper bound and descent property of function $\phi(w,w^*)$.
\begin{lemma}[Boundness of $\phi(w,w^*)$ and descent inequality]\label{Up-Gamma} Suppose Assumption~\ref{assump1} and~\ref{assump2}  hold, then the following statements also hold.
\begin{itemize}
\item[{\rm(i)}] $\phi(w,w^*)\geq d_1\|w-w^*\|^2$.
\item[{\rm(ii)}] $\phi(w,w^*)\leq d_2\|w-w^*\|^2+\epsilon [L(u,p^*)-L(u^*,p^*)]$.
\item[{\rm(iii)}] $\phi(w^k,w^*)-\mathbb{E}_{i(k)}\phi(w^{k+1},w^*)\geq d_4\|w^k-T(w^{k})\|^2+\frac{\epsilon}{N}[L(u^k,p^*)-L(u^*,p^*)]$.
\end{itemize}
\end{lemma}
The definition for global strong metric subregularity~\cite{DontchevRockafellar09} is as follows.
\begin{definition}[Global strong metric subregularity]
Let $\mathcal{H}(x)$ be a set-valued mapping between real spaces $\mathbf{X}$ and $\mathbf{Y}$. Then $\mathcal{H}(x)$ is called global strong metric subregular at $\bar{x}$ for $\bar{y}$ when $\bar{y}\in\mathcal{H}(\bar{x})$ if there exists positive number $\mathfrak{c}$ such that
\begin{equation}\label{SMS}
dist(x,\bar{x})\leq \mathfrak{c}dist\left(\bar{y},\mathcal{H}(x)\right),\quad\mbox{for all $x\in\mathbf{X}$}.
\end{equation}
\end{definition}
The following theorem regarding linear convergence of RPDC under global strong metric subregularity of the KKT mapping $H:\RR^n\times\RR^m\rightrightarrows\RR^n\times\RR^m$.
\begin{theorem}\label{prop:MSH-VIEB}{\bf (Global strong metric subregularity of $H(w)$ implies linear convergence of RPDC)}
Suppose Assumption \ref{assump1} and \ref{assump2} hold. For a given saddle point $w^*$, if $H(w)$ is global strong metric subregular at $w^*$ for $0$, then there exists $\alpha\in(0,1)$ such that
\begin{equation}
\mathbb{E}_{\mathcal{F}_{k+1}}\phi(w^{k+1},w^*)\leq\alpha^{k+1}\phi(w^0,w^*),\quad\forall k.
\end{equation}
\end{theorem}
Then the R-linear of the sequence $\{\mathbb{E}_{\mathcal{F}_{k}}w^k\}$ can be expressed as in the corollary.
\begin{corollary} [R-linear rate of $\{\mathbb{E}_{\mathcal{F}_{k}}w^k\}$] Suppose the assumptions of Theorem~\ref{prop:MSH-VIEB} hold and $\alpha\in(0,1)$ is constant. Then the sequence $\{\mathbb{E}_{\mathcal{F}_{k}}w^k\}$ converges to the desired saddle point $w^*$ at R-linear rate; i.e.,
$$\lim_{k\rightarrow\infty}\sup\sqrt[k]{\|\mathbb{E}_{\mathcal{F}_{k}}w^{k}-w^*\|}=\sqrt{\alpha}.$$
\end{corollary}
\section{Support vector machine and machine learning portfolio problem implementations}\label{num}
This section discusses experiments conducted using MATLAB R2020a on a personal computer with Intel Core i5-6200U CPU (2.40GHz) and 8.00 GB RAM. We also calculate optimal values for all experiments to check the suboptimality of RPDC using the commercial solver CPLEX 12.6.
\subsection{Support vector machine problem}
Consider the SVM problem,
\begin{eqnarray*}
\begin{array}{lll}
\mbox{\rm(SVM)}&\min\limits_{u\in[0,c]^n}   & \frac{1}{2} u^{\top}Q u-\mathbf{1}_n^{\top}u  \\
               &\mbox{s.t.}                 &  y^{\top}u=0
\end{array},
\end{eqnarray*}
where $u\in\mathbb{R}^{n}$ are the decision variables, and $Q\in\mathbb{R}^{n\times n}$ is a symmetric and positive-definite matrix. Let $Q=(Q_{1}^{\top},Q_{2}^{\top},\cdots,Q_{N}^{\top})^{\top}\in \RR^{n\times n}$ be an appropriate partition of matrix $Q$ and $Q_{i}$ be an $n_{i}\times n$ matrix. Then the KKT mapping for SVM is
\begin{equation*}
H(w)=\left(
    \begin{array}{l}
      Q u-\mathbf{1}_n+py+\mathcal{N}_{[0,c]^n}(u) \\
      y^{\top}u
    \end{array}
  \right),
\end{equation*}
where $w=\left(u,p\right)$. The following proposition shows that the KKT mapping for SVM is global strong metric subregular.
\begin{proposition}
Assume there exists at least one component $u_i^*$ of optimal solution $u^*$ that satisfies $0<u_i^*<c$. Then the KKT mapping for SVM is global strong metric subregular.
\end{proposition}
The RPDC scheme with $K(u)=\frac{1}{2}\|u\|^2$ for SVM is
\begin{eqnarray*}
&&\mbox{Choose $i(k)$ from $\{1,2,\ldots,N\}$ with equal probability}\nonumber\\
&&u^{k+1}\leftarrow\min_{u\in[0,c]^n}\langle Q_{i(k)}u^{k}, u_{i(k)}\rangle-\mathbf{1}_{n_{i(k)}}^{\top}u_{i(k)}\\
&&\qquad\quad+\langle p^k+\gamma y^{\top}u^{k},(y_{i(k)})^{\top}u_{i(k)}\rangle+\frac{1}{2\epsilon}\|u-u^k\|^2;\\
&&p^{k+1}\leftarrow p^{k}+\rho y^{\top}u^{k+1}.
\end{eqnarray*}
Thus, the primal subproblem of RPDC has the closed form
$$
\left\{
\begin{array}{l}
u_{i(k)}^{k+1}=\min\bigg{\{}\max\bigg{[}0,u_{i(k)}^k-\epsilon\bigg{(}Q_{i(k)}u^{k}-\mathbf{1}_{n_{i(k)}}\\
\qquad\qquad\qquad\qquad\qquad+(p^k+\gamma y^{\top}u^{k})y_{i(k)}\bigg{)}\bigg{]},c\bigg{\}},\\
u_{j\neq i(k)}^{k+1}=u_{j\neq i(k)}^{k}.
\end{array}
\right.$$
We used two LIBSVM datasets in the experiment: {\tt heart\_scale} ($270$ data and $13$ features) and {\tt ionosphere\_scale} ($351$ data and $34$ features). $Q$ was generated using the radial basis function kernel, and we selected $c=1$.

We partitioned the variables $N=2,5,10$ blocks, for both cases. Thus $n_i=135,54,27$ for the first dataset ({\tt heart\_scale}); and $n_i=175\;\mbox{(or $176$)},70\;\mbox{(or $71$)},35\;\mbox{(or $36$)}$ for the second dataset ({\tt ionosphere\_scale}).

\begin{figure*}
\centering
\subfigure[$\|w^k-w^*\|$]{
\begin{minipage}[b]{0.31\linewidth}
\includegraphics[width=0.95\linewidth]{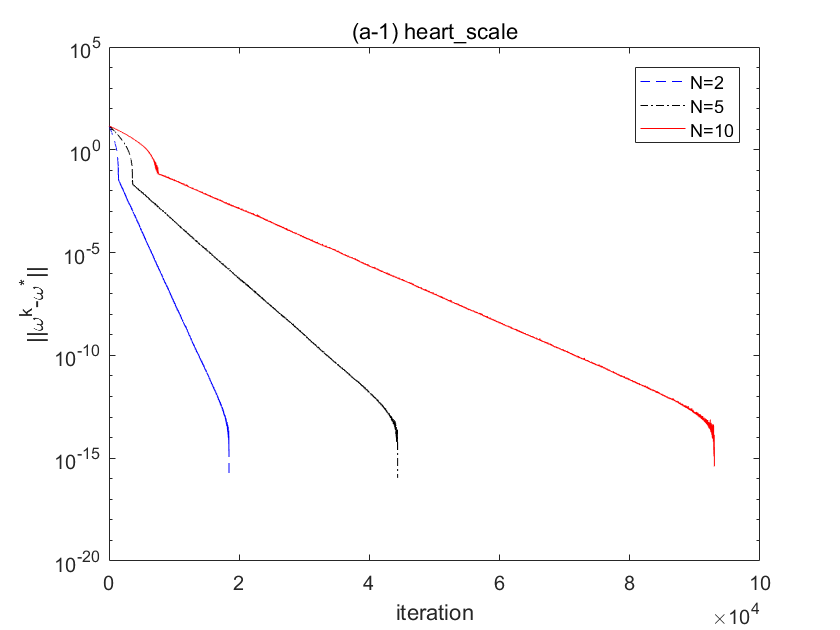}\vspace{0pt}
\includegraphics[width=0.95\linewidth]{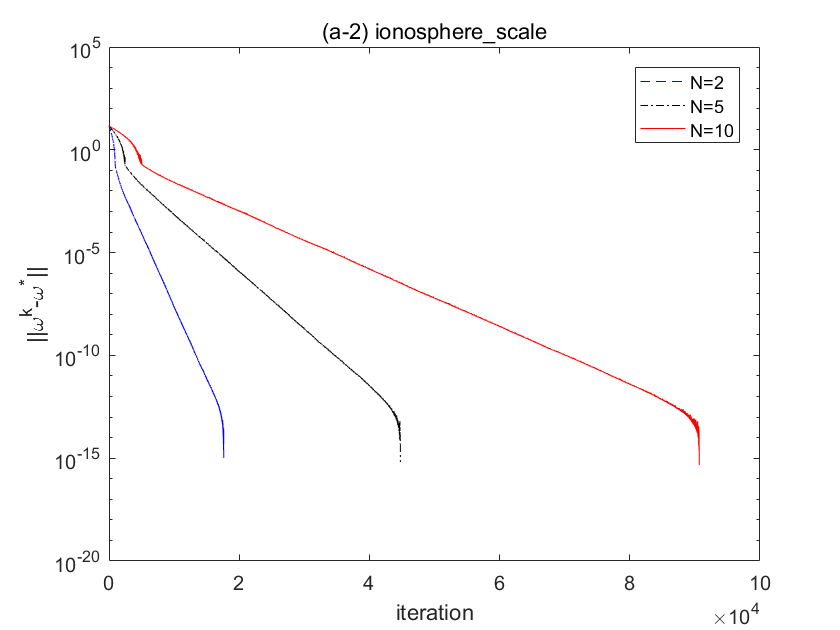}\vspace{0pt}
\includegraphics[width=0.95\linewidth]{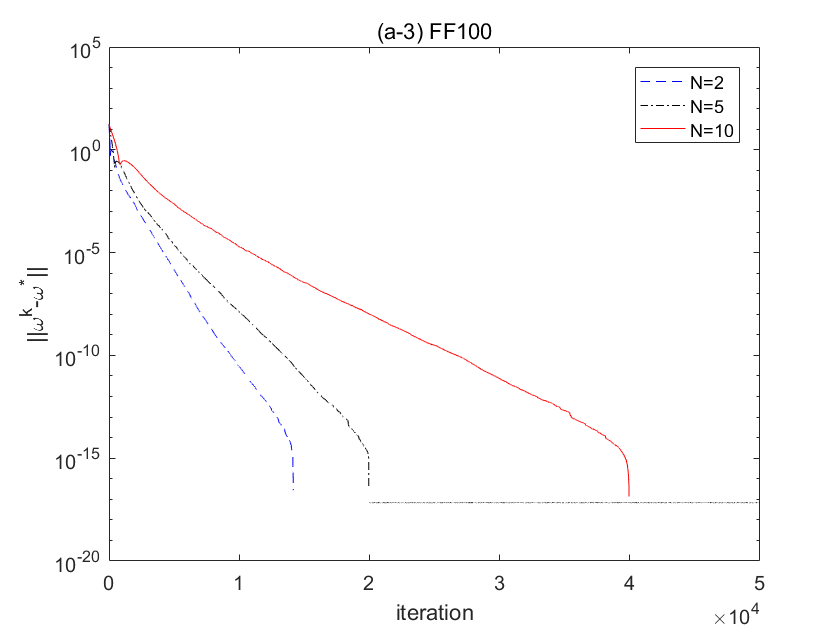}\vspace{0pt}
\includegraphics[width=0.95\linewidth]{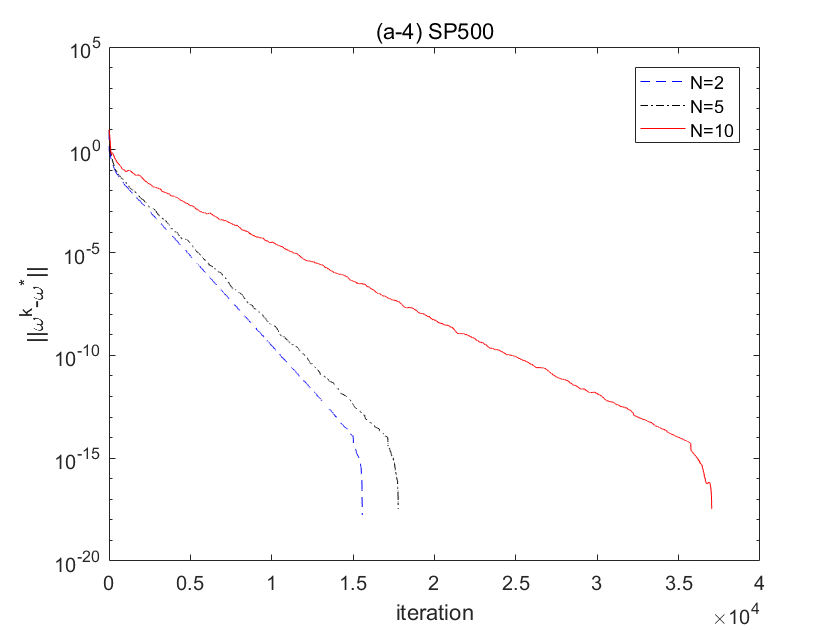}
\end{minipage}}
\subfigure[Suboptimality]{
\begin{minipage}[b]{0.31\linewidth}
\includegraphics[width=0.95\linewidth]{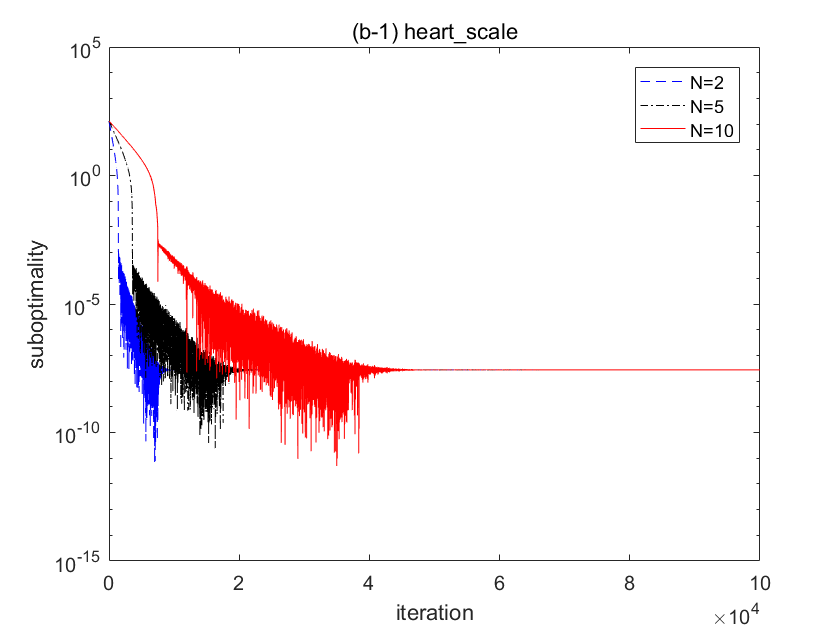}\vspace{0pt}
\includegraphics[width=0.95\linewidth]{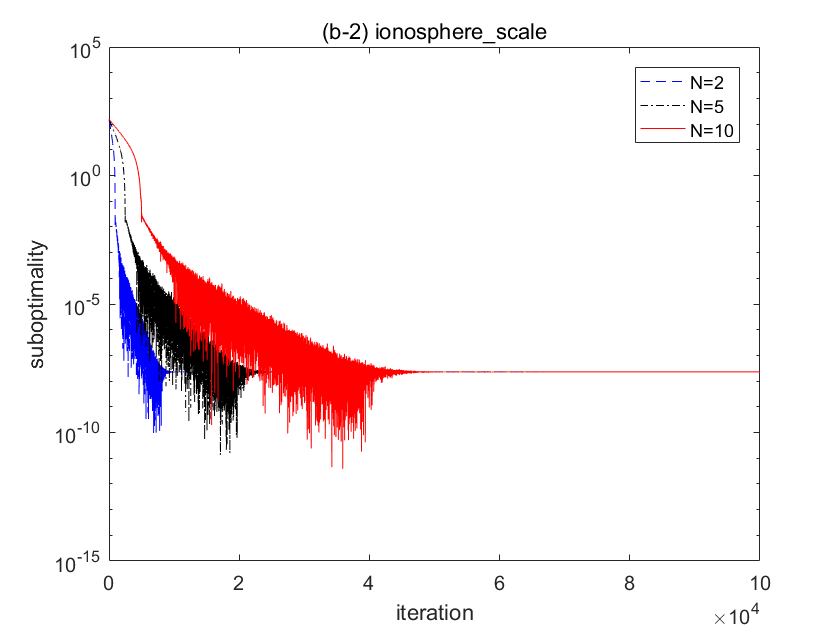}\vspace{0pt}
\includegraphics[width=0.95\linewidth]{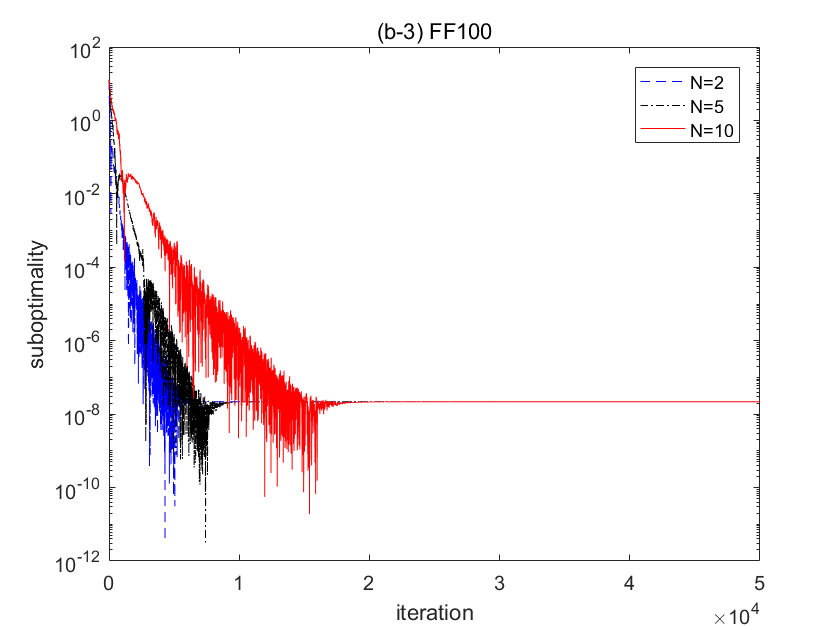}\vspace{0pt}
\includegraphics[width=0.95\linewidth]{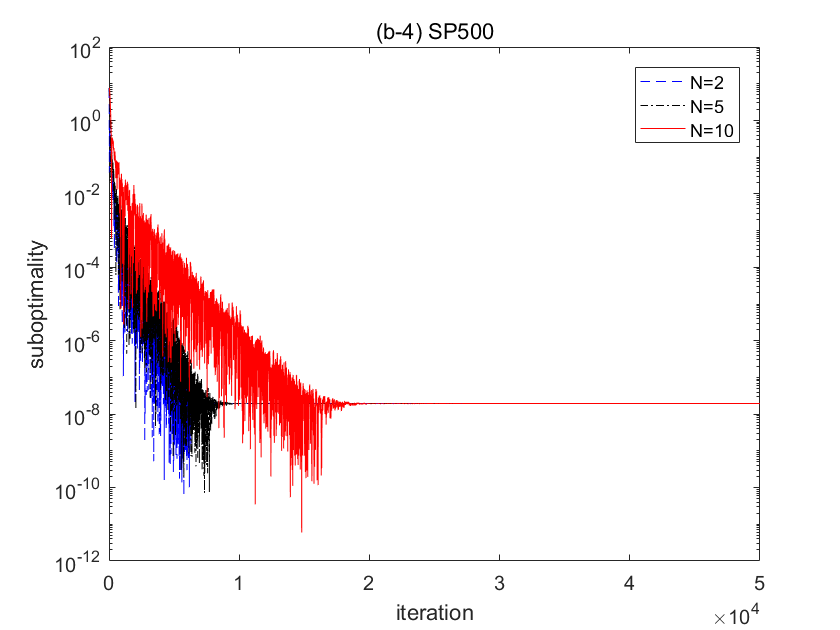}
\end{minipage}}
\subfigure[Feasibility]{
\begin{minipage}[b]{0.31\linewidth}
\includegraphics[width=0.95\linewidth]{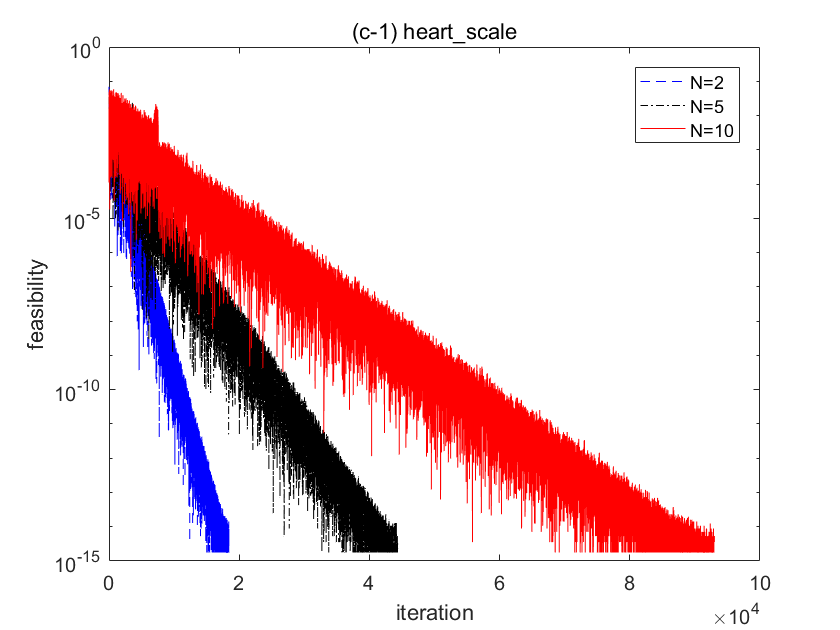}\vspace{0pt}
\includegraphics[width=0.95\linewidth]{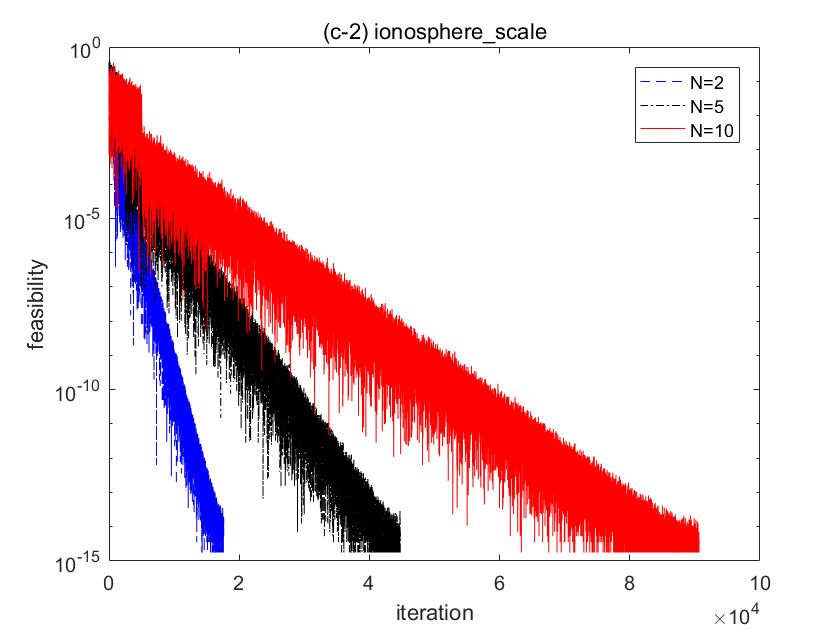}\vspace{0pt}
\includegraphics[width=0.95\linewidth]{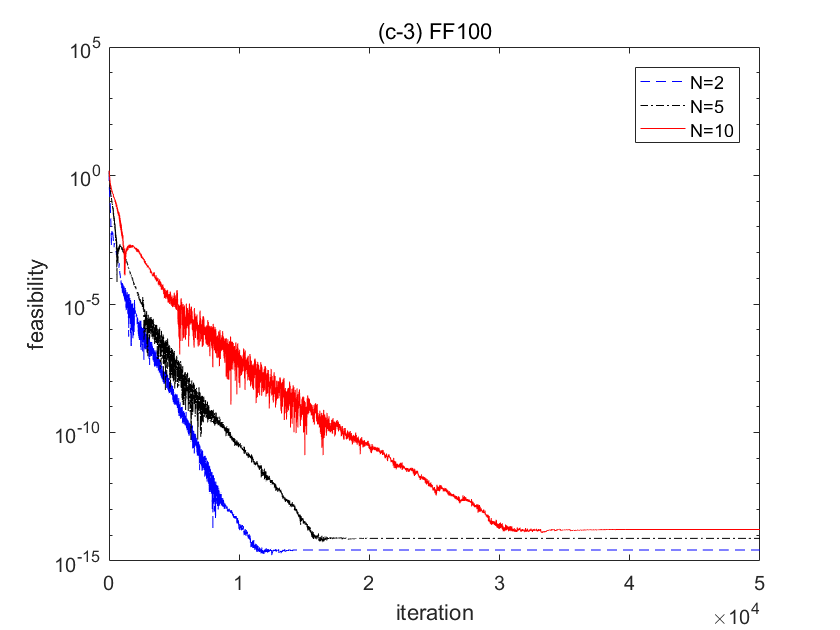}\vspace{0pt}
\includegraphics[width=0.95\linewidth]{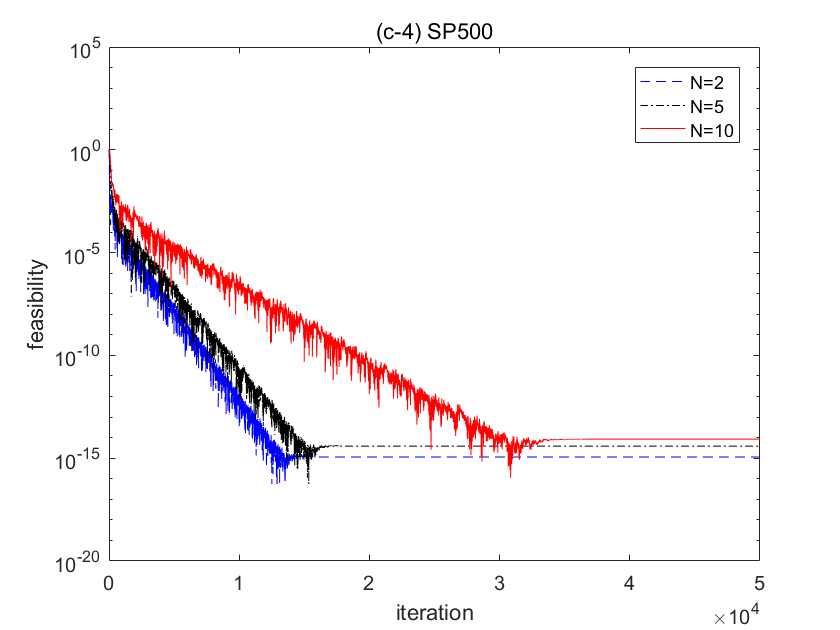}
\end{minipage}}
\caption{Number of blocks, $\|w^k-w^*\|$, suboptimality, and feasibility with respect to iteration}\label{fig:1}
\end{figure*}

In Figure~\ref{fig:1}, graphs (a-1) and (a-2) show the number of blocks and $\|w^k-w^*\|$ with respect to iteration count, respectively; graphs (b-1) and (b-2)  show the number of blocks and suboptimality with respect to iteration count, respectively; and graphs (c-1) and (c-2) show the number of blocks and feasibility with respect to iteration count, respectively.

\begin{figure*}
\centering
\subfigure[Comparison for {\tt heart\_scale}]{
\begin{minipage}[b]{0.241\linewidth}
\includegraphics[width=1\linewidth]{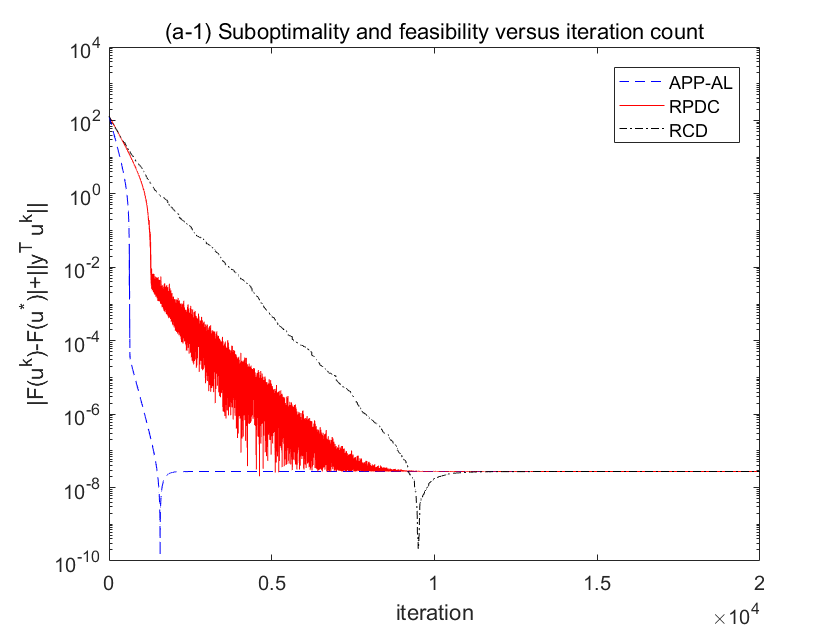}
\end{minipage}
\begin{minipage}[b]{0.241\linewidth}
\includegraphics[width=1\linewidth]{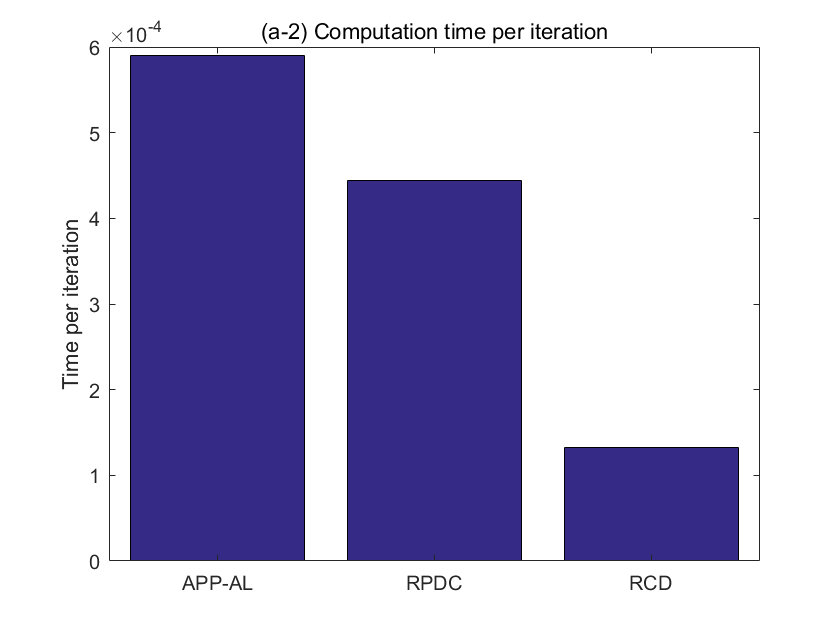}
\end{minipage}}
\subfigure[Comparison for {\tt ionosphere\_scale}]{
\begin{minipage}[b]{0.241\linewidth}
\includegraphics[width=1\linewidth]{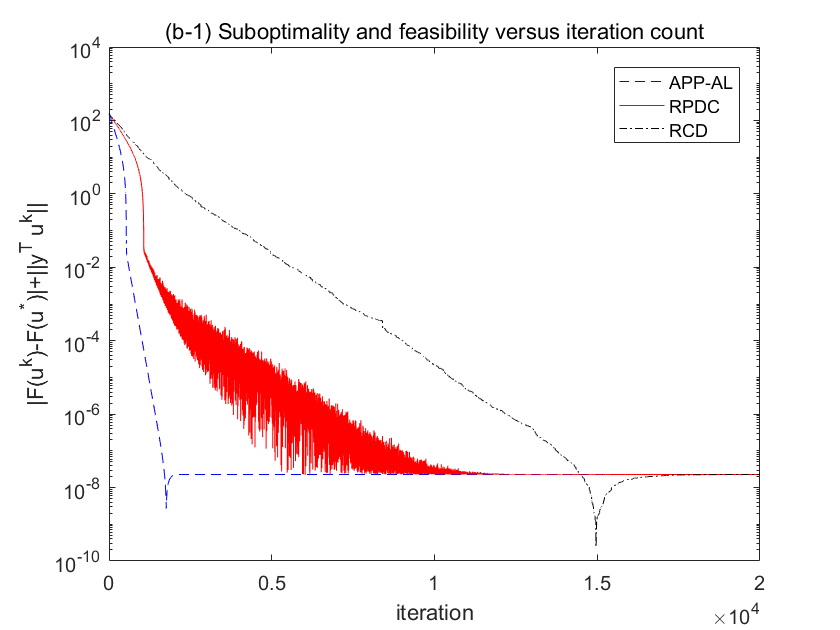}
\end{minipage}
\begin{minipage}[b]{0.241\linewidth}
\includegraphics[width=1\linewidth]{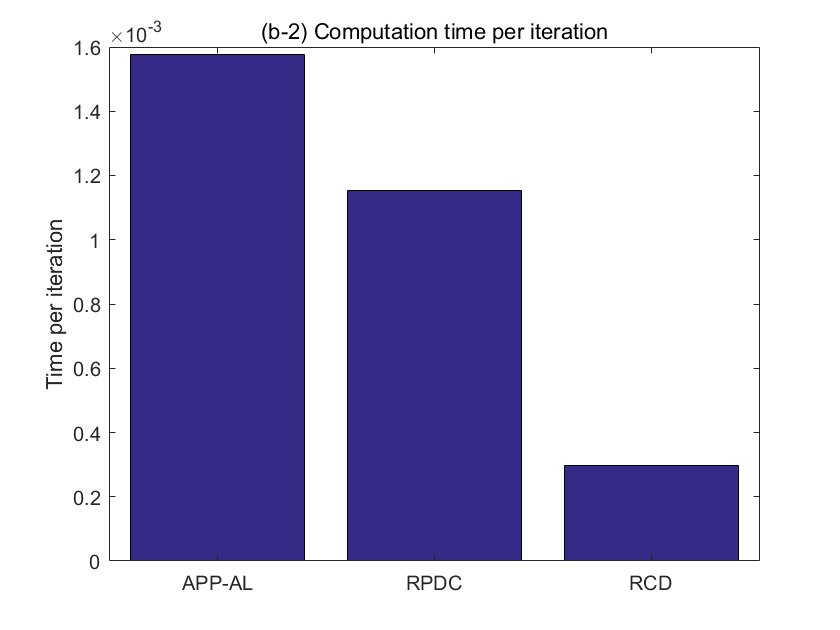}
\end{minipage}}
\caption{Comparing RPDC, APP-AL and RCD}\label{fig:2}
\end{figure*}

We compared three algorithms: APP-AL by~\cite{CohenZ} and~\cite{ZhaoZhu2017}, and RPDC from this paper with $N=2$ and random coordinate descent (RCD) algorithm~\cite{Necoara}) on {\tt heart\_scale} and {\tt ionosphere\_scale} problems. Suboptimality and feasibility were measured by $|F(u)-F(u^*)|+\|y^{\top}u\|$ with $F(u)=\frac{1}{2} u^{\top}Q u-\mathbf{1}_n^{\top}u$. In Figure~\ref{fig:2}, graphs (a-1) and (b-1) show $|F(u)-F(u^*)|+\|y^{\top}u\|$ versus iteration count; and graphs (a-2) and (b-2) show average computation time per iteration for the different algorithms. The total number of iterations required for APP-AL and RPDC are both less than RCD, APP-AL is faster than RPDC. But computation per iteration of RPDC is less than APP-AL.
\subsection{Machine learning portfolio problem}
Consider the MLP problem.
\begin{eqnarray*}
\begin{array}{lll}
\mbox{\rm(MLP)}&\min\limits_{u\in\RR^n}   & \frac{1}{2} u^{\top}\Sigma u + \lambda\|u\|_1 \\
               &\mbox{s.t.}                 &  \mu^{\top}u=\rho \\
               &                            &  \mathbf{1}_n^{\top}u=1 \\
\end{array},
\end{eqnarray*}
where $u\in \mathbb{R}^{n}$ is the decision portfolio vector, and $\Sigma\in \mathbb{R}^{n\times n}$, is the symmetric and positive-definite estimated covariance matrix of asset returns. Let $\Sigma=(\Sigma_{1}^{\top},\Sigma_{2}^{\top},\cdots,\Sigma_{N}^{\top})^{\top}\in \RR^{n\times n}$ be an appropriate partition of matrix $\Sigma$ and $\Sigma_{i}$ be an $n_{i}\times n$ matrix.

The KKT mapping for MLP is:
\begin{equation*}
H(w)=\left(
    \begin{array}{l}
      \Sigma u+\lambda\partial\|u\|_1+p_1\mathbf{1}_n+p_2\mu \\
      \mu^{\top}u-\rho\\
      \mathbf{1}_n^{\top}u-1
    \end{array}
  \right),
\end{equation*}
where $w=\left(u,p\right)$. The following proposition shows that the KKT mapping for MLP is global strong metric subregular.
\begin{proposition}
Assume there exists at least two components $u_i^*$ and $u_j^*$ for optimal solution $u^*$ that satisfy $u_i^*\neq0$ and $u_j^*\neq0$; and $\mu_i\neq\mu_j$. Then the KKT mapping for MLP is global strong metric subregular.
\end{proposition}
Therefore, the RPDC scheme with $K(u)=\frac{1}{2}\|u\|^2$ for MLP is
\begin{eqnarray*}
&&\mbox{Choose $i(k)$ from $\{1,2,\ldots,N\}$ with equal probability}\nonumber\\
&&u^{k+1}\leftarrow\min_{u\in\RR^n}\langle\Sigma_{i(k)}u^{k}, u_{i(k)}\rangle + \lambda\|u_{i(k)}\|_1\\
&&\qquad\quad+\langle p^k+\gamma\Theta(u^k),\Theta_{i(k)}(u_{i(k)})\rangle+\frac{1}{2\epsilon}\|u-u^k\|^2;\\
&&p^{k+1}\leftarrow p^{k}+\rho\Theta(u^{k+1}),
\end{eqnarray*}
where $\Theta(u)=\left(\begin{array}{l}\mu^{\top}u-\rho\\ \mathbf{1}_n^{\top}u-1\end{array}\right)$ and $\Theta_{i(k)}(u_{i(k)})=\left(\begin{array}{l}\mu_{i(k)}^{\top}u_{i(k)}\\ \mathbf{1}_{n_{i(k)}}^{\top}u_{i(k)}\end{array}\right)$.
Thus, the primal subproblem of RPDC has the closed form
$$
\left\{
\begin{array}{l}
u_{i(k)}^{k+1}=sign(\zeta_{i(k)}^k)\odot\max\{0,|\zeta_{i(k)}^k|-\epsilon\lambda\mathbf{1}_{n_{i(k)}}\},\\
u_{j\neq i(k)}^{k+1}=u_{j\neq i(k)}^{k},
\end{array}
\right.$$
where\\
$\zeta_{i(k)}^k=u_{i(k)}^k-\epsilon\left[\Sigma_{i(k)}u^{k}+\left(\mu_{i(k)},\mathbf{1}_{n_{i(k)}}\right)(p^k+\gamma\Theta(u^{k}))\right]$.

Two datasets were chosen to validate RPDC performance.
\begin{enumerate}
\item The FF100 dataset from Fama and French benchmark datasets~\cite{FF100}, created for different financial segments based on data sampled from the U.S. stock market. FF100 formed on the basis of size and book-to-market ratio; and
\item The Standard \& Poor's, USA SP500 dataset, November 2004 to April 2016, containing 442 assets and 595 observations.
\end{enumerate}
We partitioned the variables into $N=2,5,10$ blocks for both cases, i.e., $n_i=50,20,10$ and $n_i=221, 88\;\mbox{(or $90$)},\;\mbox{44 (or $46$)}$ respectively.\\
In Figure~\ref{fig:1}, graphs (a-3) and (a-4) show the number of blocks and $\|w^k-w^*\|$ with respect to iteration count; graphs (b-3) and (b-4) show the number of blocks and suboptimality with respect to iteration; and graphs (c-3) and (c-4) show the number of blocks and feasibility with respect to iteration.
\section{Conclusions}\label{ccl}
This paper proposed a randomized primal-dual coordinate (RPDC) method, a randomized coordinate extension of the first-order primal-dual method proposed by~\cite{CohenZ} and~\cite{ZhaoZhu2017}, to solve LCCP. We established almost surely convergence and expected $O(1/t)$ convergence rate for the general convex case, and expected linear convergence under global strong metric subregularity. We showed that SVM and MLP problems satisfy global strong metric subregularity under some reasonable conditions, discussed the implementation details of RPDC, and presented numerical experiments on SVM and MLP problems to verify linear convergence. Future study will consider RPDC for nonlinear convex cone programming with separable and non-separable objective and constraints.
\section*{Acknowledgements}
This research was supported by NSFC grant 71471112 and 71871140.
\nocite{langley00}


\onecolumn
\icmltitle{Supplementary material for the paper: "Linear Convergence of Randomized Primal-Dual Coordinate Method for Large-scale Linear Constrained Convex Programming"}







\icmlkeywords{Linear Constrained Convex Programming, Support Vector Machine, Machine Learning Portfolio, Randomized Coordinate Primal-dual Method, Linear Convergence}

\vskip 0.3in




First of all, we have the following observations:

In algorithm RPDC, the indices $i(k)$, $k=0,1,2,\ldots$ are random variables. After $k$ iterations, RPDC method generates a random output $(u^{k+1}, p^{k+1})$. Recall the definition of filtration $\mathcal{F}_k$ which is generated by the random variable $i(0),i(1),\ldots,i(k)$, i.e.,
$$\mathcal{F}_{k}\overset{def}{=}\{i(0),i(1),\ldots,i(k)\}, \mathcal{F}_{k}\subset\mathcal{F}_{k+1}.$$
Additionally, $\mathcal{F}=(\mathcal{F}_{k})_{k\in\mathbb{N}}$,  $\mathbb{E}_{\mathcal{F}_{k+1}}=\mathbb{E}(\cdot|\mathcal{F}_{k})$ is the conditional expectation w.r.t. $\mathcal{F}_{k}$ and the conditional expectation in term of $i(k)$ given $i(0),i(1),\ldots,i(k-1)$ as $\mathbb{E}_{i(k)}$.

Knowing $\mathcal{F}_{k-1}=\{i(0),i(1),\ldots,i(k-1)\}$, we have:
\begin{eqnarray}
\mathbb{E}_{i(k)}\langle\nabla_{i(k)}G(u^{k}),(u^{k}-u)_{i(k)}\rangle=\frac{1}{N}\langle\nabla G(u^{k}),u^{k}-u\rangle\geq\frac{1}{N}\big{[}G(u^k)-G(u)\big{]},\label{expectationG}
\end{eqnarray}
\begin{eqnarray}
\mathbb{E}_{i(k)}\big{[}J_{i(k)}(u_{i(k)}^{k})-J_{i(k)}(u_{i(k)})\big{]}=\frac{1}{N}\big{[}J(u^{k})-J(u)\big{]},\label{expectationJ}
\end{eqnarray}
and
\begin{eqnarray}
\mathbb{E}_{i(k)}\langle q^k,A_{i(k)}(u^{k}-u)_{i(k)}\rangle=\frac{1}{N}\langle q^k,A(u^{k}-u)\rangle.\label{expectationP_1}
\end{eqnarray}
Secondly, reconsidering the point $T(w^k)=\big{(}T_u(w^k),T_p(w^k)\big{)}$ generated by one deterministic iteration of APP-AL~\cite{CohenZ} for given $w^k$,
\begin{equation*}
\begin{array}{l}
\mbox{{\bf APP-AL}}\\
\left\{
\begin{array}{l}
T_u(w^k)=\arg\min\limits_{u\in \mathbf{U}}\langle\nabla G(u^{k}), u\rangle+ J(u)+\langle q^k,Au\rangle+\frac{1}{\epsilon}D(u,u^k);\\
T_p(w^k)= p^k+\gamma\left[AT_u(w^k)-b\right],
\end{array}
\right.
\end{array}
\end{equation*}
with $q^k=p^k+\gamma(Au^k-b)$, we have the following observations. The convex combination of $u^{k}$ and $T_u(w^k)$ provides the expected value of $u^{k+1}$ as following.
\begin{equation}\label{eq:u-tu-1}
\mathbb{E}_{i(k)}u^{k+1}=\frac{1}{N}T_u(w^k)+(1-\frac{1}{N})u^k,
\end{equation}
or
\begin{equation}\label{eq:u-tu}
T_u(w^k)=N\mathbb{E}_{i(k)}u^{k+1}-(N-1)u^k.
\end{equation}
Moreover, the point $T(w^k)$ satisfies that: for any $(u,p)\in\mathbf{U}\times\RR^m$,
\begin{eqnarray}\label{eq:VI-VAPP}
\left\{
\begin{array}{l}
\langle\nabla G(u^{k}),u-T_u(w^k)\rangle+J(u)-J(T_u(w^k))+\langle q^k, A(u-T_u(w^k))\rangle\\
\qquad\qquad\qquad\qquad\qquad+\frac{1}{\epsilon}\langle \nabla K(T_u(w^k))-\nabla K(u^k), u-T_u(w^k)\rangle\geq 0,\\
\gamma\left[AT_u(w^k)-b\right]=T_p(w^k)-p^k.
\end{array}
\right.
\end{eqnarray}
\newpage
\section{Proof of Lemma 1}
\proof
\quad Take $w'=w^*$ in (9), we have that
\begin{eqnarray}\label{Lambda}
\Lambda(w,w^*)&=&D(u^*,u)+\frac{\epsilon}{2N\rho}\|p-p^*\|^2+\frac{\epsilon(N-1)}{N}[L(u,p)-L(u^*,p^*)]+\frac{\epsilon(N-2)\gamma}{2N}\|Au-b\|^2\nonumber\\
&=&D(u^*,u)+\frac{\epsilon}{2N\rho}\|p-p^*\|^2+\frac{\epsilon(N-1)}{N}[L(u,p^*)-L(u^*,p^*)]+\frac{\epsilon(N-1)}{N}\langle p-p^*,Au-b\rangle\nonumber\\
&&+\frac{\epsilon(N-2)\gamma}{2N}\|Au-b\|^2.
\end{eqnarray}
\begin{itemize}
\item[{\rm(i)}]
Since $L(u,p^*)-L(u^*,p^*)\geq0$ and $\frac{1}{2\gamma}\|p-p^*\|^2+\frac{\gamma}{2}\|Au-b\|^2+\langle p-p^*,Au-b\rangle\geq0$,~\eqref{Lambda} follows that
\begin{eqnarray*}
\Lambda(w,w^*)\geq D(u^*,u)+\frac{\epsilon}{2N\rho}\|p-p^*\|^2-\frac{\epsilon(N-1)}{2N\gamma}\|p-p^*\|^2-\frac{\epsilon\gamma}{2N}\|Au-b\|^2.
\end{eqnarray*}
From Assumption 2, we have $D(u^*,u)\geq\frac{\beta}{2}\|u-u^*\|^2$. Together with the fact $Au^*=b$ and $\rho<\frac{2\gamma}{2N-1}$, above inequality follows that
$$\Lambda(w,w^*)\geq d_1\|w-w^*\|^2,$$
with $d_1=\min\bigg{\{}\frac{1}{2N}[N\beta-\epsilon\gamma\lambda_{\max}(A^{\top}A)],\frac{\epsilon}{4N\gamma}\bigg{\}}$.
\item[{\rm(ii)}]
By Young's inequality,~\eqref{Lambda} follows that
\begin{eqnarray*}
\Lambda(w,w^*)&\leq&D(u^*,u)+\frac{\epsilon}{2N\rho}\|p-p^*\|^2+\frac{\epsilon(N-1)}{N}[L(u,p^*)-L(u^*,p^*)]\nonumber\\
&&+\frac{\epsilon(N-1)}{N}[\frac{1}{2\gamma}\|p-p^*\|^2+\frac{\gamma}{2}\|Au-b\|^2]+\frac{\epsilon(N-2)\gamma}{2N}\|Au-b\|^2.
\end{eqnarray*}
From Assumption 2, we have $D(u^*,u)\leq\frac{B}{2}\|u-u^*\|^2$. Together with the fact $Au^*=b$ and $2\gamma>(2N-1)\rho$, above inequality follows that
\begin{eqnarray*}
\Lambda(w,w^*)
\leq d_2\|w-w^*\|^2+\frac{\epsilon(N-1)}{N}[L(u,p^*)-L(u^*,p^*)],
\end{eqnarray*}
with $d_2=\max\bigg{\{}\frac{(4N-3)\epsilon}{(4N-2)N\rho},\frac{NB+\epsilon(2N-3)\gamma\lambda_{\max}(A^{\top}A)}{2N}\bigg{\}}$.
\item[{\rm(iii)}] By the definition of $\Lambda(w,w')$, we have
\begin{eqnarray}
\Lambda(w,w')&\geq&\frac{\epsilon(N-1)}{N}[L(u,p)-L(u^*,p^*)]+\frac{\epsilon(N-2)\gamma}{2N}\|Au-b\|^2\nonumber\\
&=&\frac{\epsilon(N-1)}{N}[L(u,p)-L(u,p^*)]+\frac{\epsilon(N-1)}{N}[L(u,p^*)-L(u^*,p^*)]+\frac{\epsilon(N-2)\gamma}{2N}\|Au-b\|^2\nonumber\\
&\geq&\frac{\epsilon(N-1)}{N}[L(u,p)-L(u,p^*)]+\frac{\epsilon(N-2)\gamma}{2N}\|Au-b\|^2\nonumber\\
&=&\frac{\epsilon(N-1)}{N}\langle p-p^*,Au-b\rangle+\frac{\epsilon(N-2)\gamma}{2N}\|Au-b\|^2\nonumber\\
&\geq&-d_3\|p-p^*\|^2,
\end{eqnarray}
with $d_3=\frac{\epsilon(N-1)^2}{2\gamma N(N-2)}$.
\end{itemize}
\endproof
\newpage
\section{Proof of Lemma 2}
\proof\quad
{\sl Step 1: Estimate $\frac{\epsilon}{N}\mathbb{E}_{i(k)}\big{[}L(u^{k+1},q^k)-L(u,q^k)\big{]}$;}\\
For all $u\in\mathbf{U}$, the unique solution $u^{k+1}$ of the primal problem of RPDC is characterized by the following variational inequality:
\begin{eqnarray*}\label{eq:primal_VI}
\langle\nabla_{i(k)}G(u^{k}),(u^{k+1}-u)_{i(k)}\rangle+J_{i(k)}(u_{i(k)}^{k+1})-J_{i(k)}(u_{i(k)})+\langle q^k, A_{i(k)}(u^{k+1}-u)_{i(k)}\rangle\nonumber\\
+\frac{1}{\epsilon}\langle\nabla K(u^{k+1})-\nabla K(u^k), u^{k+1}-u\rangle\leq 0,
\end{eqnarray*}
which follows that
\begin{eqnarray}\label{eq:primal_VI2}
\langle\nabla_{i(k)}G(u^{k}),\big{(}u^{k}-u-(u^{k}-u^{k+1})\big{)}_{i(k)}\rangle+J_{i(k)}(u_{i(k)}^{k})-J_{i(k)}(u_{i(k)})-\big{(}J_{i(k)}(u_{i(k)}^{k})-J_{i(k)}(u_{i(k)}^{k+1})\big{)}\nonumber\\
+\langle q^k,A_{i(k)}\big{(}u^{k}-u-(u^{k}-u^{k+1})\big{)}_{i(k)}\rangle+\frac{1}{\epsilon}\langle\nabla K(u^{k+1})-\nabla K(u^k), u^{k+1}-u\rangle\leq 0.
\end{eqnarray}
Observing that for any separable mapping $\psi(u)=\sum\limits_{i=1}^{N}\psi_i(u_i)$, we have $\psi_{i(k)}(u_{i(k)}^{k})-\psi_{i(k)}(u_{i(k)}^{k+1})=\psi(u^{k})-\psi(u^{k+1})$. Therefore,~\eqref{eq:primal_VI2} follows that
\begin{eqnarray}\label{eq:primal_bound0}
&&\langle\nabla_{i(k)}G(u^{k}),(u^{k}-u)_{i(k)}\rangle+J_{i(k)}(u_{i(k)}^{k})-J_{i(k)}(u_{i(k)})+\langle q^k,A_{i(k)}(u^{k}-u)_{i(k)}\rangle\nonumber\\
&\leq&\langle\nabla G(u^{k}),u^{k}-u^{k+1}\rangle+J(u^{k})-J(u^{k+1})+\langle q^k,A(u^{k}-u^{k+1})\rangle\nonumber\\
&&+\frac{1}{\epsilon}\langle\nabla K(u^{k+1})-\nabla K(u^k), u-u^{k+1}\rangle.
\end{eqnarray}
Taking expectation with respect to $i(k)$ on both side of~\eqref{eq:primal_bound0}, together the condition expectation~\eqref{expectationG}-\eqref{expectationP_1}, we get
\begin{eqnarray}\label{eq:primal_bound1}
\frac{1}{N}\big{[}L(u^{k},q^k)-L(u,q^k)\big{]}&\leq&\mathbb{E}_{i(k)}\bigg{\{}\langle\nabla G(u^{k}),u^{k}-u^{k+1}\rangle+J(u^{k})-J(u^{k+1})\nonumber\\
&&+\langle q^k,A(u^{k}-u^{k+1})\rangle+\frac{1}{\epsilon}\langle\nabla K(u^{k+1})-\nabla K(u^k), u-u^{k+1}\rangle\bigg{\}}.
\end{eqnarray}
or
\begin{eqnarray}\label{eq:primal_bound1-1}
\frac{1}{N}\mathbb{E}_{i(k)}\big{[}L(u^{k+1},q^k)-L(u,q^k)\big{]}
&\leq&\mathbb{E}_{i(k)}\bigg{\{}\underbrace{\langle\nabla G(u^{k}),u^{k}-u^{k+1}\rangle}_{\mathfrak{a}_1}+J(u^{k})-J(u^{k+1})\nonumber\\
&&+\langle q^k,A(u^{k}-u^{k+1})\rangle+\frac{1}{N}\big{[}L(u^{k+1},q^k)-L(u^{k},q^k)\big{]}\nonumber\\
&&+\underbrace{\frac{1}{\epsilon}\langle\nabla K(u^{k+1})-\nabla K(u^k), u-u^{k+1}\rangle}_{\mathfrak{a}_2}\bigg{\}}.
\end{eqnarray}
By the gradient Lipschitz of $G$, term $\mathfrak{a}_1$ in~\eqref{eq:primal_bound1-1} is bounded by
\begin{eqnarray}\label{eq:primal_bound2}
\mathfrak{a}_1=\langle\nabla G(u^{k}),u^{k}-u^{k+1}\rangle\leq G(u^k)-G(u^{k+1})+\frac{B_G}{2}\|u^k-u^{k+1}\|^2.
\end{eqnarray}
The simple algebraic operation and Assumption 2 follows that
\begin{eqnarray}\label{eq:primal_bound4}
\mathfrak{a}_2=\frac{1}{\epsilon}\langle \nabla K(u^{k+1})-\nabla K(u^k),u-u^{k+1}\rangle
&=&\frac{1}{\epsilon}\big{[}D(u,u^k)-D(u,u^{k+1})-D(u^{k+1},u^k)\big{]}\nonumber\\
&\leq&\frac{1}{\epsilon}\big{[}D(u,u^k)-D(u,u^{k+1})\big{]}-\frac{\beta}{2\epsilon}\|u^k-u^{k+1}\|^2.
\end{eqnarray}
Combining~\eqref{eq:primal_bound1-1}-\eqref{eq:primal_bound4}, we obtain that
\begin{eqnarray}\label{eq:primal_bound5}
\frac{\epsilon}{N}\mathbb{E}_{i(k)}\big{[}L(u^{k+1},q^k)-L(u,q^k)\big{]}
&\leq&\big{[}D(u,u^k)-\mathbb{E}_{i(k)}D(u,u^{k+1})\big{]}+\mathbb{E}_{i(k)}\bigg{\{}\frac{\epsilon(N-1)}{N}\underbrace{\big{[}L(u^k,q^k)-L(u^{k+1},q^k)\big{]}}_{\mathfrak{a}_3}\nonumber\\
&&-\frac{\beta-\epsilon B_G}{2}\|u^k-u^{k+1}\|^2\bigg{\}}
\end{eqnarray}
Since $p^{k+1}=p^k+\rho(Au^{k+1}-b)$ and $q^k=p^k+\gamma(Au^{k}-b)$, term $\mathfrak{a}_3$ in~\eqref{eq:primal_bound5} follows that
\begin{eqnarray}\label{eq:primal_bound6}
\mathfrak{a}_3&=&L(u^k,q^k)-L(u^{k+1},q^k)\nonumber\\
&=&L(u^k,p^k)-L(u^{k+1},p^{k+1})+\langle q^{k}-p^{k},Au^{k}-b\rangle+\langle p^{k+1}-q^k,Au^{k+1}-b\rangle\nonumber\\
&=&L(u^k,p^k)-L(u^{k+1},p^{k+1})+\gamma\|Au^{k}-b\|^2+\rho\|Au^{k+1}-b\|^2-\gamma\langle Au^k-b,Au^{k+1}-b\rangle\nonumber\\
&=&L(u^k,p^k)-L(u^{k+1},p^{k+1})+\frac{\gamma}{2}\|Au^{k}-b\|^2+(\rho-\frac{\gamma}{2})\|Au^{k+1}-b\|^2+\frac{\gamma}{2}\|A(u^k-u^{k+1})\|^2\nonumber\\
&\leq&L(u^k,p^k)-L(u^{k+1},p^{k+1})+\frac{\gamma}{2}\|Au^{k}-b\|^2+(\rho-\frac{\gamma}{2})\|Au^{k+1}-b\|^2\nonumber\\
&&+\frac{\gamma\lambda_{\max}(A^{\top}A)}{2}\|u^k-u^{k+1}\|^2.
\end{eqnarray}
Combining~\eqref{eq:primal_bound5}-\eqref{eq:primal_bound6}, we have that
\begin{eqnarray}\label{eq:primal_bound8}
\frac{\epsilon}{N}\mathbb{E}_{i(k)}\big{[}L(u^{k+1},q^k)-L(u,q^k)\big{]}&\leq&\big{[}D(u,u^k)-\mathbb{E}_{i(k)}D(u,u^{k+1})\big{]}+\mathbb{E}_{i(k)}\bigg{\{}\frac{\epsilon(N-1)}{N}\big{[}L(u^k,p^k)-L(u^{k+1},p^{k+1})\big{]}\nonumber\\
&&-\frac{\beta-\epsilon[B_G+\frac{N-1}{N}\gamma\lambda_{\max}(A^{\top}A)]}{2}\|u^k-u^{k+1}\|^2+\frac{\epsilon\gamma(N-1)}{2N}\|Au^{k}-b\|^2\nonumber\\
&&+\frac{\epsilon(2\rho-\gamma)(N-1)}{2N}\|Au^{k+1}-b\|^2\bigg{\}}
\end{eqnarray}

{\sl Step 2: Estimate $\frac{\epsilon}{N}\mathbb{E}_{i(k)}\big{[}L(u^{k+1},p)-L(u^{k+1},q^k)\big{]}$}\\
\begin{eqnarray}\label{eq:dual_bound2-1}
L(u^{k+1},p)-L(u^{k+1},q^k)&=&\langle p-q^k,Au^{k+1}-b\rangle\nonumber\\
&=&\frac{1}{\rho}\langle p-p^{k},p^{k+1}-p^{k}\rangle-\gamma\langle Au^k-b, Au^{k+1}-b\rangle\nonumber\\
&=&\frac{1}{2\rho}\left[\|p-p^k\|^2-\|p-p^{k+1}\|^2+\|p^k-p^{k+1}\|^2\right]-\gamma\langle Au^k-b,Au^{k+1}-b\rangle\nonumber\\
&=&\frac{1}{2\rho}\left[\|p-p^k\|^2-\|p-p^{k+1}\|^2+\|p^k-p^{k+1}\|^2\right]+\frac{\gamma}{2}\|A(u^k-u^{k+1})\|^2\nonumber\\
&&-\frac{\gamma}{2}\|Au^k-b\|^2-\frac{\gamma}{2}\|Au^{k+1}-b\|^2\nonumber\\
&=&\frac{1}{2\rho}\left[\|p-p^k\|^2-\|p-p^{k+1}\|^2\right]+\frac{\gamma}{2}\|A(u^k-u^{k+1})\|^2\nonumber\\
&&-\frac{\gamma}{2}\|Au^k-b\|^2+\frac{\rho-\gamma}{2}\|Au^{k+1}-b\|^2\qquad\quad\mbox{(since $p^{k+1}=p^k+\rho(Au^{k+1}-b)$.)}\nonumber\\
&\leq&\frac{1}{2\rho}\left[\|p-p^k\|^2-\|p-p^{k+1}\|^2\right]+\frac{\gamma\lambda_{\max}(A^{\top}A)}{2}\|u^k-u^{k+1}\|^2\nonumber\\
&&-\frac{\gamma}{2}\|Au^k-b\|^2+\frac{\rho-\gamma}{2}\|Au^{k+1}-b\|^2
\end{eqnarray}
Multiply $\frac{\epsilon}{N}$ on both side of above inequality, we obtain that: $\forall p\in\mathbf{R}^m$
\begin{eqnarray}\label{eq:dual_bound4}
\frac{\epsilon}{N}\big{[}L(u^{k+1},p)-L(u^{k+1},q^k)\big{]}&\leq&\frac{\epsilon}{2N\rho}\big{[}\|p-p^k\|^2-\|p-p^{k+1}\|^2\big{]}+\frac{\epsilon\frac{1}{N}\gamma\lambda_{\max}(A^{\top}A)}{2}\|u^k-u^{k+1}\|^2\nonumber\\
&&-\frac{\epsilon\gamma}{2N}\|Au^k-b\|^2+\frac{\epsilon(\rho-\gamma)}{2N}\|Au^{k+1}-b\|^2.
\end{eqnarray}
Taking expectation with respect to $i(k)$ on both side of inequality~\eqref{eq:dual_bound4}, we have
\begin{eqnarray}\label{eq:step2}
\frac{\epsilon}{N}\mathbb{E}_{i(k)}\big{[}L(u^{k+1},p)-L(u^{k+1},q^k)\big{]}&\leq&\frac{\epsilon}{2N\rho}\big{[}\|p-p^k\|^2-\mathbb{E}_{i(k)}\|p-p^{k+1}\|^2\big{]}+\frac{\epsilon\frac{1}{N}\gamma\lambda_{\max}(A^{\top}A)}{2}\mathbb{E}_{i(k)}\|u^k-u^{k+1}\|^2\nonumber\\
&&-\frac{\epsilon\gamma}{2N}\|Au^k-b\|^2+\frac{\epsilon(\rho-\gamma)}{2N}\mathbb{E}_{i(k)}\|Au^{k+1}-b\|^2.
\end{eqnarray}
{\sl Step 3: Estimate the variance of $\Lambda(w^k,w)$.}\\
Summing inequalities~\eqref{eq:primal_bound8} and~\eqref{eq:step2}, with $d_4=\frac{\min\bigg{\{}\frac{\beta-\epsilon[B_G+\gamma\lambda_{\max}(A^{\top}A)]}{2},\frac{\epsilon[2\gamma-(2N-1)\rho]}{2N}\bigg{\}}}{\max\{N^2+2\gamma^2(N^2+2)\lambda_{\max}(A^{\top}A),4\gamma^2\}}$, we have that
\begin{eqnarray}\label{eq:sum_1}
&&\Lambda(w^k,w)-\mathbb{E}_{i(k)}\Lambda(w^{k+1},w)\nonumber\\
&\geq&\mathbb{E}_{i(k)}\bigg{\{}\frac{\epsilon}{N}\big{[}L(u^{k+1},p)-L(u,q^{k})\big{]}+\frac{\beta-\epsilon[B_G+\gamma\lambda_{\max}(A^{\top}A)]}{2}\|u^k-u^{k+1}\|^2+\frac{\epsilon[2\gamma-(2N-1)\rho]}{2N}\|Au^{k+1}-b\|^2\bigg{\}}\nonumber\\
&\geq&\mathbb{E}_{i(k)}\bigg{\{}\frac{\epsilon}{N}\big{[}L(u^{k+1},p)-L(u,q^{k})\big{]}+d_4[\left(N^2+2\gamma^2(N^2+2)\lambda_{\max}(A^{\top}A)\right)\|u^k-u^{k+1}\|^2+4\gamma^2\|Au^{k+1}-b\|^2]\bigg{\}}\nonumber\\
&\geq&\mathbb{E}_{i(k)}\bigg{\{}\frac{\epsilon}{N}\big{[}L(u^{k+1},p)-L(u,q^{k})\big{]}+d_4\big{[}\left(1+2\gamma^2\lambda_{\max}(A^{\top}A)\right)N^2\|u^k-u^{k+1}\|^2+4\gamma^2[\|A(u^k-u^{k+1})\|^2+\|Au^{k+1}-b\|^2]\big{]}\bigg{\}}\nonumber\\
&\geq&\mathbb{E}_{i(k)}\bigg{\{}\frac{\epsilon}{N}\big{[}L(u^{k+1},p)-L(u,q^{k})\big{]}+d_4\big{[}\left(1+2\gamma^2\lambda_{\max}(A^{\top}A)\right)N^2\|u^k-u^{k+1}\|^2+2\gamma^2\|Au^k-b\|^2\big{]}\bigg{\}}\nonumber\\
&=&\frac{\epsilon}{N}\mathbb{E}_{i(k)}\big{[}L(u^{k+1},p)-L(u,q^{k})\big{]}+d_4\big{[}\left(1+2\gamma^2\lambda_{\max}(A^{\top}A)\right)N^2\mathbb{E}_{i(k)}\|u^k-u^{k+1}\|^2+2\gamma^2\|Au^k-b\|^2\big{]}.
\end{eqnarray}
By Jensen's inequality,~\eqref{eq:sum_1} follows that
\begin{eqnarray}\label{eq:sum_2}
\Lambda(w^k,w)-\mathbb{E}_{i(k)}\Lambda(w^{k+1},w)&\geq&\frac{\epsilon}{N}\mathbb{E}_{i(k)}\big{[}L(u^{k+1},p)-L(u,q^{k})\big{]}\nonumber\\
&&+d_4\big{[}\left(1+2\gamma^2\lambda_{\max}(A^{\top}A)\right)N^2\|u^k-\mathbb{E}_{i(k)}u^{k+1}\|^2+2\gamma^2\|Au^k-b\|^2\big{]}.
\end{eqnarray}
Since $\mathbb{E}_{i(k)}u^{k+1}-u^k=\frac{1}{N}[T_u(w^k)-u^k]$ in~\eqref{eq:u-tu-1},~\eqref{eq:sum_2} yields that
\begin{eqnarray}\label{eq:sum_3}
\Lambda(w^k,w)-\mathbb{E}_{i(k)}\Lambda(w^{k+1},w)&\geq&\frac{\epsilon}{N}\mathbb{E}_{i(k)}\big{[}L(u^{k+1},p)-L(u,q^{k})\big{]}\nonumber\\
&&+d_4\big{[}\left(1+2\gamma^2\lambda_{\max}(A^{\top}A)\right)\|u^k-T_u(w^{k})\|^2+2\gamma^2\|Au^k-b\|^2\big{]}.
\end{eqnarray}
Since $\lambda_{\max}(A^{\top}A)\|u^k-T_u(w^{k})\|^2\geq\|A[u^k-T_u(w^{k})]\|^2$ and $T_p(w^{k})-p^k=\gamma[AT_u(w^{k})-b]$,~\eqref{eq:sum_3} follows that
\begin{eqnarray}
\Lambda(w^k,w)-\mathbb{E}_{i(k)}\Lambda(w^{k+1},w)&\geq&\frac{\epsilon}{N}\mathbb{E}_{i(k)}\big{[}L(u^{k+1},p)-L(u,q^{k})\big{]}\nonumber\\
&&+d_4\big{[}\|u^k-T_u(w^{k})\|^2+2\gamma^2\|A[u^k-T_u(w^{k})]\|^2+2\gamma^2\|Au^k-b\|^2\big{]}\nonumber\\
&\geq&\frac{\epsilon}{N}\mathbb{E}_{i(k)}\big{[}L(u^{k+1},p)-L(u,q^{k})\big{]}+d_4\big{[}\|u^k-T_u(w^{k})\|^2+\gamma^2\|AT_u(w^{k})-b\|^2\big{]}\nonumber\\
&\geq&\frac{\epsilon}{N}\mathbb{E}_{i(k)}\big{[}L(u^{k+1},p)-L(u,q^{k})\big{]}+d_4\|w^k-T(w^k)\|^2.\nonumber
\end{eqnarray}
Then we have the result of Lemma 2.
\endproof
\newpage
\section{Proof of Theorem 1 (Almost surely convergence)}
\proof
\begin{itemize}
\item[(i)] Take $w=w^*$ in Lemma 2, we have
\begin{eqnarray}\label{eq:bound2_1_7}
\Lambda(w^k,w^*)\geq\mathbb{E}_{i(k)}\Lambda(w^{k+1},w^{*})+\frac{\epsilon}{N}\mathbb{E}_{i(k)}\big{[}L(u^{k+1},p^*)-L(u^*,q^{k})\big{]}+d_4\|w^k-T(w^k)\|^2.
\end{eqnarray}
Observe that $L(u^{k+1},p^*)-L(u^*,q^{k})\geq0$. From statement (i) of Lemma 1, we have that $\Lambda(w^k,w^*)$ is nonnegative. By the Robbins-Siegmund Lemma~\cite{RS}, we obtain that $\lim\limits_{k\rightarrow+\infty}\Lambda(w^k,w^*)$ almost surely exists, $\sum\limits_{k=0}\limits^{+\infty}\|w^k-T(w^k)\|^2<+\infty$ $\mbox{a.s.}$.
\item[(ii)] Since $\lim\limits_{k\rightarrow+\infty}\Lambda(w^k,w^{*})$ almost surely exists, thus $\Lambda(w^k,w^{*})$ is almost surely bounded.
Thanks statement (i) of Lemma 1, it implies the sequences $\{w^k\}$ is almost surely bounded.
\item[(iii)] From statement (i) we have that
\begin{eqnarray*}\label{eq:ukukp1}
\lim\limits_{k\to \infty}\|w^k-T(w^k)\|=0\quad\mbox{a.s.}.
\end{eqnarray*}
By variational inequality system~\eqref{eq:VI-VAPP}, we have that any cluster point of a realization sequence generated by RPDC almost surely is a saddle point of Lagrangian for (P).
\end{itemize}
\endproof
\section{Proof of Theorem 2 (Expected   primal   suboptimality   and   expected feasibility)}
\proof
\begin{itemize}
\item[(i)] Let $h(w,w')=\Lambda(w,w')+\frac{d_3}{d_1}\Lambda(w,w^*)$. By statement (i) and (iii) in Lemma 1, we have $h(w,w')\geq0$. From Lemma 2, we obtain that
\begin{eqnarray*}\label{eq:rate_1}
\mathbb{E}_{i(k)}\frac{\epsilon}{N}\big{[}L(u^{k+1},p)-L(u,q^{k})\big{]}\leq\Lambda(w^k,w)-\mathbb{E}_{i(k)}\Lambda(w^{k+1},w)
\end{eqnarray*}
Taking expectation with respect to $\mathcal{F}_t$, $t>k$ for above inequality, we obtain that
\begin{eqnarray}\label{eq:rate_2}
\frac{\epsilon}{N}\mathbb{E}_{\mathcal{F}_t}\big{[}L(u^{k+1},p)-L(u,q^{k})\big{]}\leq\mathbb{E}_{\mathcal{F}_t}[\Lambda(w^k,w)-\Lambda(w^{k+1},w)].
\end{eqnarray}
Take $w=w^*$ in~\eqref{eq:rate_2}, we obtain
\begin{eqnarray}\label{eq:rate_2+}
0\leq\mathbb{E}_{\mathcal{F}_t}[\Lambda(w^k,w^*)-\Lambda(w^{k+1},w^*)].
\end{eqnarray}
By the combination of~\eqref{eq:rate_2} and~\eqref{eq:rate_2+}, it follows
\begin{eqnarray}\label{eq:rate_2++}
\frac{\epsilon}{N}\mathbb{E}_{\mathcal{F}_t}\big{[}L(u^{k+1},p)-L(u,q^{k})\big{]}\leq\mathbb{E}_{\mathcal{F}_t}[h(w^k,w)-h(w^{k+1},w)]
\end{eqnarray}
From the definition of $\bar{u}_t$ and $\bar{p}_t$, we have $\bar{u}_t\in\mathbf{U}$ and $\bar{p}_t\in\mathbf{R}^m$. From the convexity of set $\mathbf{U}$, $\mathbf{R}^m$ and the function $L(u',p)-L(u,p')$ is convex in $u'$ and linear in $p'$, for all $u\in\mathbf{U}$ and $p\in\mathbf{R}^m$, we have that
\begin{eqnarray}\label{eq:rate_3}
\mathbb{E}_{\mathcal{F}_t}\big{[}L(\bar{u}_{t},p)-L(u,\bar{p}_{t})\big{]}\leq\mathbb{E}_{\mathcal{F}_t}\frac{1}{t+1}\sum_{k=0}^t\big{[}L(u^{k+1},p)-L(u,q^k)\big{]}\leq\frac{Nh(w^0,w)}{\epsilon(t+1)}.
\end{eqnarray}
\item[{\rm(ii)}] If $\mathbb{E}_{\mathcal{F}_{t}}\|A\bar{u}_t-b\|=0$, statement (ii) is obviously. Otherwise, $\mathbb{E}_{\mathcal{F}_{t}}\|A\bar{u}_t-b\|\neq0$ i.e., there is set $\mathbb{W}$ such that $\mathbb{P}\{\omega\in\mathbb{W}|\|A\bar{u}_t-b\|\neq0\}>0$. Let $\hat{p}$ be a random vector:
\begin{eqnarray}\label{eq:p_omega}
\hat{p}(\omega)=\left\{
\begin{array}{cc}
0 & \omega\notin\mathbb{W} \\
\frac{M\big{(}A\bar{u}_t-b\big{)}}{\|A\bar{u}_t-b\|} & \omega\in\mathbb{W}.
\end{array}
\right.
\end{eqnarray}
Noted that for $\omega\notin\mathbb{W}$, we have $\hat{p}(\omega)=0$ and $\|A\bar{u}_t-b\|=0$. Thus
\begin{eqnarray}\label{star0}
\langle\hat{p}(\omega),A\bar{u}_t-b\rangle=M\|A\bar{u}_t-b\|=0.
\end{eqnarray}
Otherwise, for $\omega\in\mathbb{W}$, we have that
\begin{eqnarray}\label{star1}
\langle \hat{p}(\omega),A\bar{u}_t-b\rangle=M\|A\bar{u}_t-b\|.
\end{eqnarray}
Together~\eqref{star0} and~\eqref{star1}, we have
\begin{eqnarray}\label{star}
\langle\hat{p},A\bar{u}_t-b\rangle=M\|A\bar{u}_t-b\|
\end{eqnarray}
Moreover, since $Au^*=b$, we have
\begin{eqnarray}\label{eq:rate1}
L(\bar{u}_t,\hat{p})-L(u^*,\bar{p}_t)=F(\bar{u}_t)+\langle\hat{p},A\bar{u}_t-b\rangle-F(u^*)=F(\bar{u}_t)-F(u^*)+M\|A\bar{u}_t-b\|.
\end{eqnarray}
Moreover, by taking $u=\bar{u}_t$ in the right hand side of saddle point inequality, we have
\begin{eqnarray}\label{eq:rate_primal_8}
F(\bar{u}_t)-F(u^*)\geq-\langle p^*,A\bar{u}_t-b\rangle\geq-\|p^*\|\|A\bar{u}_t-b\|.
\end{eqnarray}
Combine~\eqref{eq:rate1} and~\eqref{eq:rate_primal_8}, we have that $$\|A\bar{u}_t-b\|\leq\frac{L(\bar{u}_t,\hat{p})-L(u^*,\bar{p}_t)}{\left(M-\|p^*\|\right)}.$$
Take expectation on both side of above inequality, we have that
\begin{eqnarray}
\mathbb{E}_{\mathcal{F}_{t}}\|A\bar{u}_t-b\|\leq\frac{\mathbb{E}_{\mathcal{F}_{t}}[L(\bar{u}_t,\hat{p})-L(u^*,\bar{p}_t)]}{\left(M-\|p^*\|\right)}&\leq&\mathbb{E}_{\mathcal{F}_{t}}\frac{Nh(w^0,(u^*,\hat{p}))}{\left(M-\|p^*\|\right)\epsilon(t+1)}\qquad\qquad\mbox{(by (i))}\nonumber\\
&\leq&\mathbb{E}_{\mathcal{F}_{t}}\frac{Nd_5}{\left(M-\|p^*\|\right)\epsilon(t+1)}
\end{eqnarray}
where $d_5=\sup\limits_{\|p\|<M}h(w^0,(u^*,p))$.
\item[{\rm(iii)}] Again from~\eqref{eq:rate1},~\eqref{eq:rate_primal_8} and statement (ii), statement (iii) is coming.
\end{itemize}
\endproof
\section{Proof of Lemma 3}
\proof
\begin{itemize}
\item[{\rm(i)}] This statement directly follows from the definition of $\phi(w,w^*)$ and statement (i) in Lemma 1.
\item[{\rm(ii)}] This statement directly follows from the definition of $\phi(w,w^*)$ and statement (ii) in Lemma 1.
\item[{\rm(iii)}] By the definition of $\phi(w,w^*)$, we have that.
\begin{eqnarray*}\label{lemma-bound}
&&\phi(w^k,w^*)-\mathbb{E}_{i(k)}\phi(w^{k+1},w^*)\nonumber\\
&=&\Lambda(w^k,w^*)-\mathbb{E}_{i(k)}\bigg{\{}\Lambda(w^{k+1},w^*)+\frac{\epsilon}{N}[L(u^k,p^*)-L(u^*,p^*)]-\frac{\epsilon}{N}[L(u^{k+1},p^*)-L(u^*,p^*)]\bigg{\}}\nonumber\\
&\geq&\Lambda(w^k,w^*)-\mathbb{E}_{i(k)}\bigg{\{}\Lambda(w^{k+1},w^*)+\frac{\epsilon}{N}[L(u^k,p^*)-L(u^*,p^*)]-\frac{\epsilon}{N}[L(u^{k+1},p^*)-L(u^*,q^k)]\bigg{\}}\nonumber\\
&&\qquad\qquad\qquad\qquad\qquad\qquad\qquad\qquad\qquad\qquad\qquad\qquad\qquad\qquad\qquad\mbox{(by the definition of saddle point.)}\nonumber\\
&\geq&d_4[\|w^k-T(w^{k})\|^2+\frac{\epsilon}{N}[L(u^k,p^*)-L(u^*,p^*)].\qquad\qquad\qquad\qquad\qquad\qquad\qquad\qquad\qquad\mbox{(by Lemma 2)}
\end{eqnarray*}
\end{itemize}
\endproof
\section{Proof of Theorem 3 (Global strong metric subregularity of $H(w)$ implies linear convergence of RPDC)}
\proof\quad
Considering the reference point $T(w^{k})$ associated with given point $w^k$, we have that
\begin{equation}
\left\{
\begin{array}{l}
0\in\nabla G(u^k)+\partial J(T_u(w^k)) +A^{\top} q^k+\frac{1}{\epsilon}\left[\nabla K(T_u(w^k))-\nabla K(u^{k})\right]+\mathcal{N}_{\mathbf{U}}(T_u(w^k))      \\
0=b-AT_u(w^k)+\frac{1}{\gamma}\left[T_p(w^k)-p^{k}\right]
\end{array}
\right.
\end{equation}
Thus
\begin{eqnarray*}
\begin{array}{l}
v(T(w^k))=\left(
\begin{array}{l}
\nabla G(T_u(w^k))-\nabla G(u^k) +A^{\top}(T_p(w^k)-q^k)+\frac{1}{\epsilon}\left[\nabla K(u^{k})-\nabla K(T_u(w^k))\right]      \\
\frac{1}{\gamma}\left[p^k-T_p(w^k)\right]
\end{array}
\right)\in H(T(w^{k})).
\end{array}
\end{eqnarray*}
From Assumption 1 and 2, there is $\delta>0$ such that
\begin{eqnarray}\label{MS-VEB}
\|v(T(w^k))\|^2\leq\delta\|w^k-T(w^k)\|^2.
\end{eqnarray}
Since $H(w)$ is global strong metric subregular at $w^*$ for $0$, then
\begin{eqnarray}\label{MS+}
\|T(w^{k})-w^*\|\leq\mathfrak{c}dist(0,H(T(w^{k})))\leq \mathfrak{c}\|v(T(w^k))\|\leq\mathfrak{c}\sqrt{\delta}\|w^k-T(w^{k})\|.
\end{eqnarray}
Since $\|w^{k}-w^*\|\leq\|T(w^{k})-w^*\|+\|w^{k}-T(w^{k})\|$, we have
\begin{equation}\label{Strong-VAPP-EB}
\|w^{k}-w^*\|\leq (\mathfrak{c}\sqrt{\delta}+1)\|w^k-T(w^{k})\|.
\end{equation}
From statement (iii) of Lemma 3, we have that
\begin{eqnarray}\label{linear_1}
\phi(w^{k},w^*)-\mathbb{E}_{i(k)}\phi(w^{k+1},w^*)&\geq&d_4\|w^k-T(w^{k})\|^2+\frac{\epsilon}{N}[L(u^k,p^*)-L(u^*,p^*)]\nonumber\\
&\geq&\frac{d_4}{(\mathfrak{c}\sqrt{\delta}+1)^2}\|w^{k}-w^*\|^2+\frac{\epsilon}{N}[L(u^k,p^*)-L(u^*,p^*)]\qquad\mbox{(by~\eqref{Strong-VAPP-EB})}\nonumber\\
&\geq&\delta'\{d_2\|w^{k}-w^*\|^2+\epsilon [L(u^k,p^*)-L(u^*,p^*)]\}\nonumber\\
&\geq&\delta'\phi(w^{k},w^*).\qquad\qquad\qquad\qquad\qquad\qquad\qquad\mbox{(by~(i) of Lemma 3)}
\end{eqnarray}
where $\delta'=\min \{\frac{d_4}{\max\{d_2(\mathfrak{c}\sqrt{\delta}+1)^2,d_4+1\}}, \frac{1}{N+1}\}<1$. It follows that
\begin{equation}
\mathbb{E}_{i(k)}\phi(w^{k+1},w^*)\leq\alpha\phi(w^k,w^*).
\end{equation}
where $\alpha=1-\delta'\in(0,1)$. Taking expectation with respect to $\mathcal{F}_{k+1}$ for above
inequality, we obtain that
\begin{equation}
\mathbb{E}_{\mathcal{F}_{k+1}}\phi(w^{k+1},w^*)\leq\alpha^{k+1}\phi(w^0,w^*).
\end{equation}
\endproof
\section{Proof of Corollary 1 (R-linear rate of the sequence $\{\mathbb{E}_{\mathcal{F}_{k}}w^k\}$)}
\proof
\quad By statement (i) in Lemma 3, we have that $\phi(w,w^*)\geq d_1\|w-w^*\|^2$. By Theorem 3, we have that
\begin{equation*}
\mathbb{E}_{\mathcal{F}_{k}}\phi(w^{k},w^*)\leq\alpha^{k}\phi(w^0,w^*).
\end{equation*}
Then we have that
\begin{eqnarray*}
\mathbb{E}_{\mathcal{F}_{k}}\|w^k-w^*\|^2\leq\frac{\alpha^{k}\phi(w^0,w^*)}{d_1}.
\end{eqnarray*}
By convexity of $\|\cdot\|^2$ and Jensen's inequality, we obtain that
\begin{eqnarray*}
\|\mathbb{E}_{\mathcal{F}_{k}}w^k-w^*\|\leq\hat{M}(\sqrt{\alpha})^{k}\quad\mbox{with}\;\hat{M}=\sqrt{\frac{\phi(w^0,w^*)}{d_1}}.
\end{eqnarray*}
This shows that the sequence $\{\mathbb{E}_{\mathcal{F}_{k}}w^k\}$ converges to the desired saddle point $w^*$ at R-linear rate; i.e.,
$$\lim_{k\rightarrow\infty}\sup\sqrt[k]{\|\mathbb{E}_{\mathcal{F}_{k}}w^{k}-w^*\|}=\sqrt{\alpha}<1.$$
\endproof
\newpage
\section{Proof of Proposition 1}
\proof
\quad By the piecewise linear of $H(w)$ and Zheng and Ng~\cite{Zheng2014}, we have that $H(w)$ is global metric subregular at $w^*$ for $0$. Since $Q$ is positive-definite, then problem (SVM) has unique solution $u^*$. Hence, to show $H(w)$ is global strongly metric subregular, we need to prove uniqueness of the Lagrangian multiplier for (SVM). Suppose their are two multipliers $p$ and $p'$, thus we have
\begin{eqnarray*}
\left\{\begin{array}{l}
0\in Q u^*-\mathbf{1}_n+ {p}y+\mathcal{N}_{[0,c]^n}(u^*)  \\
0\in Q u^*-\mathbf{1}_n+ {p}'y+\mathcal{N}_{[0,c]^n}(u^*)
\end{array}
\right.
\end{eqnarray*}
Since there exists at least one component $u_i^*$ of optimal solution $u^*$ satisfies $0<u_i^*<c$, then $\xi_i =\mathcal{N}_{[0,c]}(u_i^*)=0$. Thus, we have that
\begin{eqnarray}
\left\{
\begin{array}{l}
Q_i u^*-1+y_ip =0\\
Q_i u^*-1+y_ip' =0
\end{array}
\right.
\end{eqnarray}
We conclude that $p = p'$. Therefore $H(w)$ is global strongly metric subregular.
\endproof

\section{Proof of Proposition 2}
\proof
\quad By the piecewise linear of $H(w)$ and Zheng and Ng~\cite{Zheng2014}, we have that $H(w)$ is global metric subregular at $w^*$ for $0$. Since $\Sigma$ is positive-definite, then problem (MLP) has unique solution $u^*$. Hence, to show $H(w)$ is global strongly metric subregular, we need to prove uniqueness of the Lagrangian multiplier for (MLP). Suppose their are two pare of multipliers $(p_1,p_2)$ and $(p_1',p_2')$, thus we have
\begin{eqnarray*}
\left\{\begin{array}{l}
0\in \Sigma u^* +\lambda \partial \|u^* \|_1 + {p}_1\mu+{p}_2\mathbf{1}_n  \\
0\in \Sigma u^* +\lambda \partial \|u^* \|_1 + {p}_1'\mu+{p}_2'\mathbf{1}_n
\end{array}
\right.
\end{eqnarray*}
Since $u_i^*\neq0$, $u_j^*\neq0$, thus $\xi_i = \partial |u_i^*|$ and $\xi_j=\partial|u_j^*|$ are single valued and we have
\begin{eqnarray}
\left\{
\begin{array}{l}
\Sigma_i u^*+\lambda\xi_i+\mu_ip_1+p_2 =0\\
\Sigma_i u^*+\lambda\xi_i+\mu_ip_1'+p_2' =0
\end{array}
\right.\\
\left\{
\begin{array}{l}
\Sigma_j u^*+\lambda\xi_j+\mu_jp_1+p_2 =0\\
\Sigma_j u^*+\lambda\xi_j+\mu_jp_1'+p_2' =0
\end{array}
\right.
\end{eqnarray}
It follows that
\begin{eqnarray}
\left\{
\begin{array}{l}
\mu_i(p_1-p_1')+p_2-p_2' =0\\
\mu_j(p_1-p_1')+p_2-p_2' =0
\end{array}
\right.
\end{eqnarray}
Since $\mu_i\neq\mu_j$, we conclude that $p_1 = p_1'$ and $p_2=p_2'$. Therefore $H(w)$ is global strongly metric subregular.
\endproof
\nocite{langley00}
\newpage
\bibliography{paper}
\bibliographystyle{icml2020}

\end{document}